\documentstyle[11pt]{article}
\begin{document}
\title{\rm The Theorem of Kronheimer-Mrowka}
\author{\rm Vishwambhar Pati,
\\ Indian Statistical Institute\\ Bangalore 560059, India}
\date{}
\maketitle

\newcommand{\mod}[1]{\mid\!{#1}\!\mid}
\newcommand{\Lam}{\Lambda}
\newcommand{\lam}{\lambda}
\newcommand{\tet}{\theta}
\newcommand{\sym}[3]{#1^{#2_{#3}}}
\newcommand{\til}[1]{\tilde{#1}}
\newcommand{\summ}[3]{\sum_{#1=#2}^{#3}}
\newcommand{\cdel}{c~-~\frac{1}{2\pi}[(\rho_{g}^{-1}\delta)_{{\cal H}_{g}}]}
\newcommand{\inte}[2]{\int_{#1}^{#2}}
\newcommand{\kernel}[1]{\mbox{Ker}\,{#1}}
\newcommand{\coker}[1]{\mbox{Coker}\,{#1}}
\newcommand{\im}[1]{\mbox{Im}\,{#1}}
\newcommand{\norml}{\left\|}
\newcommand{\normr}{\right\|}
\newcommand{\brce}[2]{_{\{#1,#2\}}}
\newcommand{\fr}[2]{\frac{#1}{#2}}
\newcommand{\alp}{\alpha}
\newcommand{\bet}{\beta}
\newcommand{\gam}{\gamma}
\newcommand{\del}{\delta}
\newcommand{\eps}{\epsilon}
\newcommand{\sig}{\sigma}
\newcommand{\smy}[3]{#1_{#2}^{#3}}
\newtheorem{lemma}{Lemma}[subsection]
\newtheorem{proposition}[lemma]{Proposition}
\newtheorem{theorem}[lemma]{Theorem}
\newtheorem{notation}[lemma]{\bf Notation :}
\newtheorem{remark}[lemma]{Remark}
\newtheorem{corollary}[lemma]{Corollary}
\newtheorem{exercise}[lemma]{Exercise}
\newtheorem{example}[lemma]{Example}
\newcommand{\bol}[2]{{\bf #1}^{#2}}
\newcommand{\gras}[3]{{\bf #1}(#2,#3)}
\newtheorem{definition}[lemma]{Definition}
\newcommand{\rarr}{\rightarrow}
\newcommand{\larr}{\leftarrow}
\newcommand{\mtxii}[4]{\left(\begin{array}{cc}#1&#2\\#3&#4\end{array}\right)}
\def\R{{\rm I\!\!R}}
\def\K{{\rm I\!\!K}}
\def\H{{\rm I\!\!H}}
\def\Proj{{\rm I\!\!P}}
\def\C{{\rm \kern.24em \vrule width .02em
   height1.4ex depth -.05ex \kern -.26em C}}
\def\Z{{\bf Z}}
\newcommand{\rot}[1]{#1^{\frac{1}{2}}}
\newcommand{\iroot}[1]{#1^{-\frac{1}{2}}}
\newcommand{\noin}{\noindent}
\newcommand{\zc}{\bar{z}}
\newcommand{\tr}{\mbox{tr}}
\newcommand{\fun}[2]{\pi_{1}(#1,#2_{0})}
\newcommand{\pair}[3]{(#1, #2_{#3})}
\newcommand{\boun}{\partial}
\newcommand{\ext}{\bigwedge}

\section{Some Transversality Results}

\subsection{Introduction} 
This note is an exposition of the proof of Thom's
Conjecture, namely that algebraic curves minimise genus within their
homology class in $\C\Proj^{2}$, due to Kronheimer-Mrowka (see [KM]).
New invariants for 4-manifolds, called Seiberg-Witten invariants, are
used in the proof. 

For the sake of completeness we have included an appendix at the end
on Fredholm theory. As background material, the reader may wish to consult the references 
[D],[PP], [W]. 

\subsection{Metrics}
Let $X$ be a connected, compact, oriented 4-manifold. Fix a
reference metric $g_{0}$ on $X$, which defines a volume form $dV_{g_{0}}$
on $X$, and hence a trivialisation of the determinant bundle 
$\ext^{4}(T^{*}X)$ of the real cotangent bundle. 
This trivialisation of $\ext^{4}$ will be fixed throughout, and will be 
denoted $dV$ without any ambiguity.

Thus the structure group for $X$ is reduced to $SL(4,\R)$, and
one lets $P_{X}SL(4)$ denote the principal $SL(4,\R)$ bundle on
$X$ consisting of all frames in $TX$ on which $dV$ yields the
constant function 1 on $X$. 

The corresponding ad-bundle, denoted $\mbox{ad}\underline{sl}_{4}$,
the vector bundle on $X$ whose fibre over $x~\in~X$  is the
vector space of traceless endomorphisms $\mbox{End}^{0}(T_{x}X)$
of $T_{x}X$ splits into the direct sum of $\mbox{ad}\underline{so}_{4}$
and $\mbox{ad}{\cal P}$ (abuse of notation since ${\cal P}$ isn't
a lie algebra) corresponding to the Cartan decomposition :
$\underline{sl}_{4}=\underline{so}_{4}\oplus {\cal P}$. Here, the bundle 
 $\mbox{ad}\underline{so}_{4}$ is the adjoint bundle corresponding to
the principal $SO(4,\R)$-bundle $P_{X}SO(4)$ consisting of $g_{0}$-orthonormal
 oriented frames. It has fibre consisting of traceless
$g_{0}$-skew-symmetric endos of $T_{x}X$ over $x$. $\mbox{ad}{\cal P}$
has fibre consisting of traceless
$g_{0}$-symmetric endos of $T_{x}X$ over $x$. 
Both are real rank-3 bundles.

Clearly, if $g$ is {\em another} Riemannian metric on $X$ with $dV_{g}=dV$,
pointwise, then $g(u,v)=g_{0}((\exp h)u,(\exp h)v)$ for some
section $h$ of  $\mbox{ad}{\cal P}$, and thus the space of all
$C^{r}$ Riemannian metrics $g$ whose volume form coincides with
the prescribed one $dV$ is precisely the space of
$C^{r}$-sections $\Gamma^{r}(\mbox{ad}{\cal P})$. This space
clearly contains exactly one representative from each conformal
class of $C^{r}$-metric on $X$, and may consequently be thought
of as the space of equivalence classes of conformal
$C^{r}$-metrics on $X$. 

The space ${\cal C}:=\Gamma^{r}(\mbox{ad}{\cal P})$ maybe given a Banach
space structure (the reference metric $g_{0}$ is the origin, corresponds
to the zero section) via the norm

$$\norml s\normr := \sup\left\{\norml \nabla^{i}_{g_{0},v} s\normr_{x,g_{0}}:
0\;\leq i\leq\;r, \;\;v\in\otimes^{i}T\;X,\;\norml v\normr_{x,g_{0}}=1\right\}
$$
which is the topology of uniform convergence of all covariant
derivatives upto order $r$. We will fix $r$ to be suitably large
later on. 

\subsection{$\ast_{g}$ and Self-Duality}

In whatever follows, {\em all} metrics $g$ will be elements of
${\cal C}$, unless indicated otherwise.

Let $\Omega^{i}(X)$ denote the space of smooth sections of 
$\ext^{i}(T^{*}X)$. The bundle $\ext^{i}(T^{*}X)$ will be denoted
 simply as $\ext^{i}(X)$ in future.
The bundle $\ext^{4}(X)$ is identified with $\ext^{0}(X)$,
the trivial bundle on $X$, via the trivialisation defined
in the last section, i.e., the volume form $dV=dV_{g}$ goes to the constant 
function 1 under this identification. 

The Hodge star-operator $*_{g}$ is the pointwise operator which makes the 
following diagram commute:
\begin{eqnarray}
\begin{array}{ccccc}
\ext^{2}(X)&\otimes &\ext^{2}(X)&\stackrel{\wedge}{\longrightarrow}&
\ext^{0}(X)\\
\left |\!\right |&  & \uparrow\,*_{g}& &\left |\!\right |\\
\ext^{2}(X)&\otimes &\ext^{2}(X)&\stackrel{(\;,\; )_{g}}
{\longrightarrow}&\ext^{0}(X)
\end{array}
\label{stardiagram}
\end{eqnarray}
where $(\;,\;)_{g}$ denotes the pointwise $g$ inner product on 2-forms.
One easily checks that $*_{g}~\circ~*_{g}~=~\mbox{id}$, and the
fact that $*_{g}$ is an isometry with respect to the global $g$-inner
product defined on $\Omega^{2}(X)$ by :
$$\langle\omega, \tau\rangle_{g} = \int_{X} (\omega,\tau)_{g} dV =\int_{X}\omega\wedge
*_{g}\tau $$
\begin{remark}{\rm We remark here that {\em on 2-forms} the
$*$-operator is a conformal invariant of the metric. This follows
because if $g'=\lam g$, where $\lam$ is a $C^{\infty}$ function, then
$(\omega,\tau)_{g'}=\lam^{-2}(\omega,\tau)_{g}$ whereas
$dV_{g'}=\lam^{2}dV_{g}$, so that 
$(\omega,\tau)_{g'}dV_{g'}=(\omega,\tau)_{g}dV_{g}$. This
shows that the $*$-operator and the global inner product $<\,,\,>_{g}$
are both invariant under conformal changes of the metric. }
\label{conformal}
\end{remark}

We then have the eigenbundle decomposition with respect to the involution
$*_{g}$, namely 

$$\wedge^{2}(X)=\wedge^{2+}_{g} \oplus \wedge^{2-}_{g}$$
where $\ext^{2\pm}_{g}$ denotes the $\pm 1$-eigenspaces of
$*_{g}$ in $\ext^{2}(X)$. The projections $\pi^{\pm}_{g}(\omega)$ are
$\frac{1}{2}(\omega\pm *_{g}\omega)$, and $g$-orthogonal projections
with respect to the pointwise $g$-inner product. The 
space $\Omega^{2}(X)$ 
then decomposes correspondingly to $\Omega^{2\pm}_{g}(X)$,
and this decomposition is orthogonal with respect to the global 
inner-product $\langle\;,\;\rangle_{g}$. 

\begin{lemma} 
(Decomposition for $\Omega^{2+}_{g}$) {\rm There
is a $\langle\;,\;\rangle_{g}$-orthogonal decomposition :

$$ \Omega^{2+}_{g} = \mbox{Im}d_{g}^{+} \oplus {\cal H}^{2+}_{g}$$
where $d_{g}^{+}:=\pi^{+}_{g}\circ d$ and ${\cal H}^{2+}_{g}$ denotes
the $\Delta_{g}$-harmonic forms in $\Omega^{2+}_{g}$.}
\label{decomp}
\end{lemma}

\noin
{\bf Proof:} By the Hodge decomposition theorem, if 
$\alp\in\Omega^{2+}_{g}$, we may write 

$$\alp = \alp_{1} + \alp_{2} + \alp_{3} $$
where $\alp_{1}\in {\cal H}^{2}_{g}$, $\alp_{2}\in
\mbox{Im}\,d$, and
$\alp_{3}\in \mbox{Im}\,\delta_{g}$, where $\delta_{g} =
*_{g}d*_{g}$ is the adjoint of $d$ with respect to the global
$L^{2}$-inner product $<\;,\;>_{g}$, and the decomposition is 
orthogonal with respect thereto. Since $*_{g}$ and $\Delta_{g}$ commute,
it follows that $*_{g}$ maps ${\cal H}_{g}^{2}$ to itself, and interchanges
$ \mbox{Im}\,\delta_{g}$ and  $\mbox{Im}\,d_{g}$. Thus
$\alp=*_{g}\alp$ implies that $\alp_{1}=*_{g}\alp_{1}$, 
$\alp_{2}=*_{g}\alp_{3}$, and $\alp_{3}=*_{g}\alp_{2}$. Writing 
$\alp_{2}=d\gamma$, we have $\alp = \alp_{1} + d\gam + *_{g}d\gam$,
and since $d\gam + *_{g}d\gam = 2\pi^{+}_{g}d\gam =
2d^{+}_{g}\gam$, the lemma follows. \hfill $\Box$.

\vspace{.25in}
In future, it will be our convention to identify $H^{2}(X,\R)$
with the space of $g_{0}$-harmonic 2-forms 
${\cal H}^{2}_{g_{0}}$, which will henceforth be simply written 
${\cal H}^{2}$. Similarly, the symbols $*$, $\pi^{\pm}$, $\Omega^{\pm}$, 
${\cal H}^{2\pm}$, 
and $d^{\pm}$ without subscripts will mean that the reference metric
$g_{0}$ is understood. The commutative square ~(\ref{stardiagram}) above 
implies that $<\alp,\bet>_{g}=\int_{X}\alp\wedge *_{g}\bet = 
[\alp]\cup*_{g}[\bet]$ for
$\alp ,\bet\in {\cal H}^{2}_{g}$ and all Riemannian 
metrics $g$. Box brackets around a closed form $\alp$ 
will always denote cohomology class (=$\Delta_{g_{0}}$-harmonic 
component).

\begin{remark}{\rm If $\alp\in {\cal H}^{2+}_{g}$, and $\alp~\neq~0$, then 
$[\alp]\cup [\alp]~=~ <\alp,*_{g}\alp>_{g}=\norml\alp\normr_{g}^{2} ~>~ 0$, 
and similarly 
$\alp\in {\cal H}^{2-}_{g}$ implies $[\alp]\cup [\alp] ~<~0$. Thus if
$\alp$ is self (resp. anti-selfdual) with respect to {\em any}
metric $g$, its cup product with itself is positive (resp. negative).
So, for any metric $g$, the cup product pairing is positive
definite (resp. negative definite) on ${\cal H}^{2+}_{g}$ (resp.
${\cal H}^{2-}_{g}$). Hence the numbers $b_{2}^{\pm}=
\mbox{dim}{\cal H}^{2\pm}_{g}$ give the signature type of the cup product 
pairing $\cup$ for
any metric $g$. }
\label{cupsign}
\end{remark}
\vspace{.25in}
\begin{proposition}{\rm Let $\mbox{dim}H^{2}(X,\R)=r$,
$b_{2}^{+}=1$, and $b_{2}^{-}=r-1$. Let the positive cone of the
$\cup$-product pairing be denoted by
$$C:=\{\alp\in H^{2}(X,\R):\; \alp \cup\alp > 0 \}$$
and let $C_{+}$ and $C_{-}$ denote the two components of $C$ determined
by the sign of $\pi^{+}\alp$. Then for any metric $g\in {\cal
C}$, there exists a unique $g$-harmonic 2-form
$\omega_g\in {\cal H}^{2+}_{g}$ satisfying
(i) $[\omega_{g}] \in C_{+}$, (ii) $[\omega_{g}]\cup [\omega_{g}] = 
\langle\omega_{g},*_{g}\omega_{g}\rangle~=~1$. Finally, (iii) $[\omega_{g}]=
[\omega_{g'}]$
if $g$ and $g'$ are conformally equivalent.
}
\label{omegag}
\end{proposition}
\vspace{.25in}
\noindent
{\bf Proof:} For any Riemannian metric $g$, 
 there is the composite map :
$${\cal H}^{2}\hookrightarrow \Omega^{2} \rarr {\cal H}^{2}_{g}$$
which we call $\psi_{g}$. For all $g$, this is an isomorphism, taking a
cohomology class (=$\Delta_{g_{0}}$-harmonic form) to its $\Delta_{g}$-
harmonic representative. Clearly $\psi^{-1}_{g}(\alp)=[\alp]$ for
$\alp\in {\cal H}_{g}^{2}$. Via $\psi^{-1}_{g},\;*_{g}$ maybe be regarded as an
involution of $H^{2}(X,\R)={\cal H}^{2}$. 
 Choose a $<\;,\;>_{g}$ unit length element of the
$+1$-eigenspace of $*_{g}$ (which is one dimensional by assumption) 
 in $H^{2}(X,\R)$. Changing the sign of
this element if necessary, one can ensure that it lies in
$C_{+}$. This cohomology class has unique $g$-harmonic
representative, denoted by $\omega_{g}$, and is the required element, 
proving (i) and (ii). 
Finally, (iii) follows from the Remark~\ref{conformal}.\hfill $\Box$

\subsection{The Map $\rho_{g}$}
\label{rhog}
We recall the bundle isometry 
$\rho_{g}:\Lambda^{2+}_{g}\rarr \underline{su}(W_{+})$ of \S 1.4 in [PP],
where $W_{\pm}$ are the Hermitian rank 2 bundles arising from the 
$\mbox{Spin}_{c}$ structure on $X$ compatible with the Riemannian metric
$g$. The metrics on the bundles $h_{\pm}$ is to be fixed once and for all,
and thus the corresponding (Hilbert-Schmidt) inner product on the bundles  
$\underline{su}^(W_{+})$ are also fixed once and for all. Let us make 
the isometric identification 
$$\rho_{0}:\Lambda^{2+}_{g_{0}}:=\Lambda^{2+}
\rarr\underline{su}(W_{+})$$ 
where $\rho_{0}:=\rho_{g_{0}}$, 
with respect to the reference metric $g_{0}$ once and for all. Then
we can view $\rho_{g}$ as a  bundle isometry 
$$\rho_{g}:\Lambda^{2+}_{g}\;\rarr \;\Lambda^{2+}$$
with $\rho_{0}=\mbox{Id}$.
This map can be explicitly described as follows. 
Write
$$g(v,w)~=~ g_{0}\left( \exp(h)v,\,\exp(h)w \right)$$
 where $h$ is a smooth
section of the bundle $\mbox{ad}\cal{P}$ of \S 1.2 above. Then for 
$\omega,\;\tau \in \Lambda^{2}(T^{*}_{x}(X))$, one has 
$$(\omega, \tau)_{g} = 
(\wedge^{2}(\exp(-h)\omega,\wedge^{2}(\exp(-h)\tau)_{0}$$
Thus $\rho_{g}=\wedge^{2}(\exp(-h))$. 

Note that $\ext^{4}(\exp(-h))=\ext^{2}(\rho_{g})=\mbox{Id}$, 
 since $dV_{g}=dV_{g_{0}}=dV$, 
i.e. $h$ is a traceless endomorphism of the 
tangent bundle. Now one can calculate
the Hodge star operator $*_{g}$ in terms of $* :=*_{g_{0}}$.

\begin{lemma}{\rm  The Hodge star operators $*_{g}$ and $*$ are related
by the formula :
$$\rho_{g}*_{g} =  * \rho_{g}$$
Consequently, $\rho_{g} \pi^{+}_{g}=\pi^{+} \rho_{g}$, and the bundle
map $\rho_{g}:\ext^{2} \rarr \ext^{2}$ is an automorphism of the bundle
$\ext^{2}$, effecting an isometry between the 
metrics $(\;,\;)_{g}$ and $(\;,\;)_{0}$, carrying $\ext^{2+}_{g}$ onto
$\ext^{2+}$ and the $g$-self-dual 2-forms $\Omega^{2+}_{g}$ to the 
$g_{0}$-self-dual 2-forms $\Omega^{2+}$.
}
\end{lemma}
\label{starstar}

\vspace{.25in}
\noin
{\bf Proof:} 
\begin{eqnarray*}
\omega\wedge_{*_{g}}\tau = (\omega, \tau)_{g}dV_{g}&=&(\rho_{g}\omega,
\rho_{g}\tau)_{0}dV \\
&=& \rho_{g}\omega\wedge *\rho_{g}\tau =
 \wedge^{2}(\rho_{g})(\omega \wedge\rho_{g}^{-1}*\rho_{g}\tau) \\
&=& \omega\wedge\rho_{g}^{-1}*\rho_{g}\tau
\end{eqnarray*}
since $\wedge^{2}(\rho_{g})=\mbox{Id}$.
Thus $*_{g}=\rho_{g}^{-1}*\rho_{g}$, and the result follows.
\hfill$\Box$

\subsection{Transversality}
\vspace{.25in}
\begin{notation}{\rm

Let $b_{2}^{+}\geq 1$. Let ${\cal C}$ be as in the  subsection 1.2,
and, as per our convention, $\Omega^{2+},\; {\cal H}^{2+}$ be
the $+1$ eigenspaces with respect to $*_{g_{0}}$ where $g_{0}$
is the reference metric. Let $\pi_{\cal H}$ denote the orthogonal 
(with respect to $\langle \;,\;\rangle_{0}$) 
projection from $\Omega^{2+}$ to ${\cal H}^{2+}$  from Lemma 
~\ref{decomp}.  Let $\delta \in \Omega^{2+}$, and let
$c=c_{1}(L)$ denote a fixed element of $H^{2}(X,\R)=
{\cal H}^{2}\subset\Omega^{2}$ such that 
$c\cup c ~<~0$. 

Let $G$ denote the Grassmanian of $b_{2}^{-}$-dimensional
 subspaces of ${\cal H}^{2}$. There is a natural rank 
$b_{2}^{+}$ bundle $\gamma^{+}$ on $G$ whose fibre over 
$P \in G$ 
is ${\cal H}^{2}/P$.  
One may regard this as the quotient bundle of $G\times {\cal H}^{2}$
by the tautological subbundle on $G$. 
By definition, there is the natural quotient map 
$\tau :G\times{\cal H}^{2}\rarr
\gamma^{+}$ of bundles on $G$, viz. $\tau(P,\alp)=\alp \pmod{P}$. 

If $c\in {\cal H}^{2}$ is a cohomology class, then $c$ defines the 
constant section of $G\times {\cal H}^{2}$, also denoted by $c$. 
Let  its image $\tau (c)$ be denoted by $ s_{c}$, a section of 
$\gamma^{+}$. The zero locus of this section is :
$$S_{c}=\{P \in G: c\in P\}$$

If $g\in {\cal C}$ is a Riemannian metric, we have the $g$-self dual 
projector $\pi^{+}_{g}:{\cal H}^{2}_{g}~\rarr~{\cal H}^{2+}_{g}$.
 Thus one has a natural map :

\begin{eqnarray*}
P&:&{\cal C}\rarr  G\\
 & & g \mapsto \mbox{Ker}(\pi^{+}_{g}\circ\psi_{g}: {\cal H}^{2}\rarr
{\cal H}^{2+}_{g})
\end{eqnarray*}
where $\psi_{g}$ is the isomorphism from ${\cal H}^{2}$ to 
${\cal H}^{2}_{g}$ introduced in ~\ref{omegag}. This map is 
 easily seen to be $C^{r}$.

Note that if $P(g)\in S_{c}$, we have that 
$\psi_{g}(c)\in {\cal H}^{2,-}_{g}$, which implies $[\psi_{g}(c)]\cup
[\psi_{g}(c)]= c \cup c < 0$,
 by ~\ref{cupsign}. So $\im{P}\cap S_{c}= \phi$ if $c\cup c > 0$.
When $c\cup c < 0$,  
 We have the following lemma, proved in
\S 5.4 of [PP] (see also [DK], \S 4.3.14):

\vspace{.25in}
\noin 
{\bf Lemma:} Let $c\in {\cal H}$ satisfy $c \cup c < 0$. The map 
$P~:~{\cal C}~\rarr~G$ defined above is transverse to $S_{c}$, and its
 inverse image (the submanifold of `$c$-bad metrics') 
$B_{c}=P^{-1}(S_{c})$ is therefore 
 a $C^{r}$ submanifold  of ${\cal C}$ of codimension $b_{2}^{+}$. 

\vspace{.25in}
Now let $\delta\in \Omega^{2+}$. Then $\rho_{g}^{-1}(\delta)\in 
\Omega^{2+}_{g}$. Its $g$-harmonic component 
$\rho^{-1}_{g}(\delta)_{{\cal H}_{g}}$ is in ${\cal H}^{2+}_{g}$, and 
thus defines a cohomology class $[\rho^{-1}_{g}(\delta)_{{\cal H}_{g}}]$
 in $H^{2}(X,\R)={\cal H}^{2+}$. 
Consider the element $c(g,\delta)=
\cdel \in {\cal H}=H^{2}(X,\R)$. (Box brackets denote cohomology class.)
Note $c(g, 0)\equiv c$. We then have the section 
$(g, c(g,\delta))$ of the trivial bundle ${\cal C}\times {\cal H}^{2}$. 
Let $\sigma_{c,\delta}$ denote the image of this section under the bundle
map $P^{*}\tau~:~{\cal C}\times {\cal H}~\rarr~P^{*}(\gamma^{+})$. Note
that $\sigma_{c,0}$ is just the section $P^{*}s_{c}$ where 
$s_{c}$ is defined above. Let us denote $\sigma_{c,0}$ by $\sigma_{c}$.

Let $B_{c,\delta}$ denote the zero locus of this section. 
Note $B_{c,0}$ is just the submanifold $B_{c}$ defined above. 

Let ${\cal I}$ be a $C^{r}$-embedded arc,
parametrised by $[-1, 1]$ in ${\cal C}$, meeting $B_{c}$ transversely.
We will denote a typical element of ${\cal I}$ by $g_{t}$. We will also
use the subscript $t$ wherever a subscript $g_{t}$ occurs, e.g.
$\omega_{t}$ for $\omega_{g_{t}}$ etc., for notational convenience.
We now have a couple of transversality lemmas, for the two separate
cases $b_{2}^{+}\geq 2$ and $b_{2}^{+}=1$. 
}
\label{transnotation}
\end{notation}

\vspace{.25in}
\begin{proposition}{\rm 
Let $b_{2}^{+}\geq 2$, and ${\cal I}$ be transverse to $B_{c}$ as above.
In this case of $b^{+}_{2}\geq 2$, this means ${\cal I}\cap B_{c}~=~\phi$.
Then, for all $\delta$ in an $\eps$-ball $U=B(0,\eps)$ of the origin in
$\Omega^{2+}$, the intersection ${\cal I}\cap B_{c,\delta}$ is empty.
}
\label{transverse1}
\end{proposition}

\vspace{.25in}
\noin
{\bf Proof:} We are given that ${\cal I}\cap B_{c}=\phi$. This means
that ${\cal I}$ does not meet the zero locus of the section $\sigma_{c}=
\sigma_{c,0}$ defined above, i.e. that 
$\sigma_{c,0}(g_{t})\neq 0$ for all $t\in [-1, 1]$. If $\norml\;\normr$ is
some bundle metric on $P^{*}(\gamma^{+})$, let 
$a:=\mbox{min}(\norml\sigma_{c,0}(g_{t})\normr:t \in [-1,1])$, so that
$a > 0$. Since $\sigma_{c,\delta}(g_{t})$ varies continuously with 
$\delta$,  
if one chooses $\eps$ small enough, then one can arrange that
$\norml \sigma_{c,0}(g_{t})-\sigma_{c,\delta}(g_{t})\normr <
\frac{a}{2}$ for all $t\in {\cal I}$ and all $\delta\in B(0,\eps)$. 
It will then follow that
$\sigma_{c,\delta}(g_{t})$ is non-zero for all $t\in [-1,1]$, i.e. that
${\cal I}\cap B_{c,\delta} =\phi $, proving the proposition.\hfill
$\Box$.

\begin{proposition}{\rm 
Let $b_{2}^{+}=1$, and let ${\cal I}$ meet the submanifold $B_{c}$ 
transversally at a single point, say $\{g_{a}\}$. Then, for all $\delta$
in an $\eps$-ball $U=B(0,\eps)$ of the origin in $\Omega^{2+}$, the 
map 
\begin{eqnarray*}
f_{\delta}&:& {\cal I}\rarr \R\\
& & g_{t}a \mapsto \left(\cdel\right)\cup [\omega_{g}]=c(g,\delta)\cup
[\omega_{g}]
\end{eqnarray*}
where $[\omega_{g}]=\psi_{g}^{-1}(\omega_{g})$ is the cohomology class
of $\omega_{g}\in {\cal H}_{g}^{2+}$, has a unique zero, 
and this zero is a regular value. 
}
\label{derivative}
\end{proposition}

\vspace{.25in}
\noin
{\bf Proof:} Let $\sigma_{c}(g_{t})$ denote the section $\sigma_{c}$ 
restricted to the arc ${\cal I}$. Let 
$p~:~P^{*}(\gamma^{+})~\rarr~ {\cal C}$ denote the bundle projection,
and $Z$ the zero section of this bundle. Since $b_{2}^{+}=1$, this a 
real line bundle, and may be regarded as the line subbundle of 
${\cal C}\times\Omega^{2}$ whose fibre over $g$ is ${\cal H}_{g}^{2+}$.  
This bundle has a trivialisation over all of ${\cal
C}$, defined by $g\mapsto \omega_{g}$, where 
$\omega_{g}\in {\cal H}^{2+}_{g}\subset \Omega^{2+}_{g}$
is the form defined in ~\ref{omegag}.  By the fact that the zero section
of $\gamma^{+}$ and the section $s_{c}(G)$ meet transversely inside 
$\gamma^{+}$, and the transversality of $P$ to $S_{c}$, it easily
follows that $\sigma_{c}({\cal C})$ and $Z$ meet transversely in 
$P^{*}(\gamma^{+})$. By assumption, the intersection: 
$${\cal I}\cap B_{c}= p\sigma_{c}({\cal I})\cap 
p(\sigma_{c}({\cal C})\cap Z)$$ 
is a transverse intersection at the single point $\{g_{a}\}$. 
Since $p:\sigma_{c}({\cal C})\rarr {\cal C}$ is a 
diffeomorphism, which carries $\sigma_{c}({\cal I})$ and 
$\sigma_{c}({\cal C})\cap Z$ diffeomorphically to ${\cal I}$ and $B_{c}$ 
respectively, it follows that inside the manifold
$\sigma_{c}({\cal C})$, $\sigma_{c}({\cal I})$ meets 
$\sigma_{c}({\cal C})\cap Z$ transversely 
at the single point 
$\sigma_{c}(g_{a})$. This means that the  
tangent space $T_{\sigma_{c}(g_{a})}(\sigma_{c}({\cal I}))$, which is 
a one dimensional subspace of 
$T_{\sigma_{c}(g_{a})}(\sigma_{c}({\cal C}))$
 is linearly independent of  the subspace 
$T_{\sigma_{c}(g_{a})}(\sigma_{c}({\cal C})\cap Z)=
 T_{\sigma_{c}(g_{a})}(\sigma_{c}({\cal C}))~\cap T~_{\sigma_{c}(g_{a})}(Z)$, 
which is of 
codimension one in $T_{\sigma_{c}(g_{a})}(\sigma_{c}({\cal C}))$. This means that the 
tangent space $T_{\sigma_{c}(g_{a})}(\sigma_{c}({\cal I}))$ is not 
contained in $T_{\sigma_{c}(g_{a})}(Z)$. Since this last space is of
codimension one in the tangent space 
$T_{\sigma_{c}(g_{a})}(P^{*}\gamma^{+})$, one has that the curve 
$\sigma_{c}({\cal I})$ meets the zero section $Z$ transversally in 
the singleton $\{\sigma_{c}(g_{a})\}$. Now let 
$\theta~:~P^{*}(\gamma^{+})\rarr\R$ be the map 
$\theta(\alp)=[\alp]\cup [\omega_{g}]=\int_{X}\alp\wedge\omega_{g}$, 
where $\alp\in {\cal H}_{g}^{+}$, 
coming from the trivialisation of the line bundle $P^{*}(\gamma^{+})$
described above. Then $\theta$ is a submersion, and $Z=\theta^{-1}(0)$. 
So saying that the curve 
$\sigma_{c}({\cal I})$ meets the zero section $Z$ transversally in 
the singleton $\{\sigma_{c}(g_{a})\}$ is equivalent to the statement 
that $f_{0}=\theta\circ\sigma_{c}:{\cal I}\rarr \R$ has a unique zero at
$g_{a}$, and $f_{0}'(a):=f_{0}'(g_{a})$ is non-zero. Say $f_{0}'(a) >0$.
 Then,
\begin{description}
\item[(i)]$f_{0}'(t)$ will be strictly positive on a closed neighborhood 
$V$ of $g_{a}$, and also 

\item[(ii)] $\min\{|f_{0}(g_{t})|:g_{t}\in {\cal I}-V^{\circ}\}=b>0$
\end{description}

Since $\sigma_{c,\delta}$ varies smoothly with $\delta$, (i) and (ii)
above continue to hold good for 
$f_{\delta}:=\theta\circ\sigma_{c,\delta}$ for all $\delta\in U$, 
where $U$ is an $\eps$-ball around $0$ for $\eps$ suitably small.
This proves the proposition.\hfill $\Box$

\section{The Parametrised Seiberg-Witten Moduli Space}
\subsection{The Gauge Group Action}
In the sequel, we assume $b_{2}^{+}(X)\geq 1$. Let $L,\, W_{+},\, W_{-}$
be a $\mbox{Spin}_{c}$ structure on $X$. $c$ will always denote
$c_{1}(L)\in H^{2}(X,\R)={\cal H}^{2+}$. Let $g_{0}$ be the
reference metric as before. Let $k\geq 6$. We need to define various spaces :
\begin{description}
\item ${\cal A}:=$ the completion of $C^{\infty}\;\;
U(1)$-connections on $L$ (an affine space modelled on $\Omega^{1}(X)$),
with respect to Sobolev k-norm $L^{2}_{k}$. By abuse of
language, we shall denote the $L^{2}_{k}$-completion of $\Omega^{1}$
as $\Omega^{1}$. 
\item $\Gamma(W_{+}):=$ $L^{2}_{k}$ completion of
$C^{\infty}$ complex-valued sections of $W_{+}$.
\item $\Gamma(W_{-}):=$ $L^{2}_{k-1}$ completion of
$C^{\infty}$ complex-valued sections of $W_{-}$.
\item ${\cal G}:=\,\mbox{Map}(X, S^{1})$. A Hilbert manifold whose lie algebra 
is the $L^{2}_{k+1}$-completion of $\Omega^{0}(X)$.
\item $\Omega^{2+}:=\,L^{2}_{k-1}$ completion of real valued
$*_{g_{0}}$ self-dual 2-forms, again denoted by same symbol by abuse
of language. 
\item ${\cal N}:=\,{\cal A}\times \Gamma(W_{+})$.
\item ${\cal N}^{*}:=\,{\cal A}\times (\Gamma(W_{+})-\{0\})$.
\end{description}

The choice of $k\geq 6$ implies, by Sobolev's Lemma, that elements of
each of the Sobolev spaces defined above are at least twice
continuously differentiable.

${\cal G}$ acts on ${\cal N}$ by the action :
$$g.(A,\Phi)= (A - (\frac{1}{2\pi i})g^{-1}dg,\, g\Phi)= (gA, g\Phi)$$
and the choice of norms above makes this a smooth action. Note
that $g.(A,0)=(A,0)$ for all $g\in S^{1}\subset {\cal G}$, so
that the gauge group action of ${\cal G}$ on ${\cal N}$ is free
on ${\cal N}^{*}$, and has (as can be checked easily) 
 isotropy $S^{1}$ on $(A,0)$.
 
\begin{remark}{\rm The decomposition lemma ~\ref{decomp} for 
$\Omega^{2+}=\Omega^{2+}_{g_{0}}$
continues to hold for the Sobolev completions defined above
because $d:\Omega^{1}_{L^{2}_{k}}\rarr \Omega^{2}_{L^{2}_{k-1}}$
has closed range and ${\cal H}^{2+}_{L^{2}_{k-1}}$ consists of smooth
forms by elliptic regularity. }
\end{remark}
\vspace{.25in}
Let ${\cal I}$ denote a compact arc in ${\cal C}$, meeting $B_{c}$ 
as in the hypotheses of 
Propositions ~\ref{transverse1} and ~\ref{derivative}. 
A typical point on ${\cal I}$
will be denoted by $g_{t}$, and a subscript $t$ anywhere will
mean that the metric $g_{t}$ is being used. Absence of subscript, 
as always, will mean that the reference metric
$g_{0}$ is understood. We also recall the map $f_{\delta}:{\cal I}\rarr
\R$ in the setting of $b_{2}^{+}=1$, that was defined in 
Corollary~\ref{derivative}.  For simplicity let us assume that
$f_{\delta}^{-1}(0)$ is the reference metric  $g_{0}$, 
and $0$ is a regular value as stated there. 
\vspace{.25in}
We recall the notation and definitions of \S 3, 4 and 5 of [PP]. 
For a metric $g_{t}$ on $X$, $(A,\Phi)$ is called a
`monopole' if it satisfies the following equations: 
\begin{eqnarray}
D_{A,t}\Phi = 0\nonumber\\
\rho_{t}(F_{A}^{+,t})-\sigma(\Phi,\Phi)=0
\label{monopole}
\end{eqnarray}
where $F^{+,t}_{A}=\pi^{+}_{t}(F_{A})$, and 
we have, once and for all, identified 
$\Gamma(\underline{su}(W_{+}))$ 
with $\Omega_{0}^{2+}=\Omega^{2+}$ via $\rho_{0}$, (see
~\ref{rhog}), and the pairing
$\sigma : W_{+}\otimes \overline{W}_{+}\rarr \Omega^{2+}\simeq 
\Gamma(\underline{su}(W_{+}))$ is the pairing defined in [PP], \S
1.5. $\rho_{t}=\rho_{g_{t}}$ is the isomorphism identifying 
$\Omega^{2+}_{t}$ with $\Gamma(\underline{su}(W_{+}))=\Omega^{2+}$
as defined in ~\ref{rhog} (see also \S 1.4 of [PP]). 

Let ${\cal G}$ act trivially on $\Omega^{2+}$ and ${\cal I}$, and
consider the $G$-equivariant map :
\begin{eqnarray*}
Q : &{\cal A}\times \Gamma (W_{+})\times {\cal I}&\rarr \;\Gamma (W_{-})
\times \Omega^{2+}\\
&(A,\Phi, g_{t})&\mapsto \;
(D_{A,t}\Phi, \rho_{t}(F_{A}^{+,t})-\sigma (\Phi,\Phi))
\end{eqnarray*}
 
Thus the  solutions to the Seiberg-Witten equations
(\ref{monopole}) are precisely
the elements of $Q^{-1}(0,0)$. Since $(0,0)$ won't in general be
a regular value for $Q$, we need to consider the
$\delta$-perturbed Seiberg-Witten equations, viz.
$Q^{-1}(0,\delta)$, where $\delta$ is a suitably small element
of $\Omega^{2+}$. Since the map $Q$ is ${\cal G}$ equivariant,
and $(0,\delta)$ is fixed by ${\cal G}$, it is natural to
quotient out $Q^{-1}(0,\delta)$ by the ${\cal G}$-action. We
shall proceed to do this in detail, in the sequel.

\subsection{The Derivative of Q}
\begin{lemma}{\rm  The derivative of $Q$ is as follows :
\begin{description}
\item[(i)] 
$$DQ_{(A,\Phi,g_{t})}(a,\phi, 0) 
= \left(D_{A,t}\phi + 2\pi ia\circ \Phi,\; 
\rho_{t}(d_{t}^{+}a)-2\Im\sigma (\Phi,\phi)
\right)$$
where $\Im$ denotes imaginary part, and $\circ$ Clifford multiplication. 
\item[(ii)] In the case when $b_{2}^{+}=1$, let 
$\delta := Q(A,0,g_{0})=F_{A}^{+}$, and ${\cal I},\; f_{\delta}$
be as in Proposition ~\ref{derivative}. Then:
$$DQ_{A,0,g_{0}}(0,0,g'(0))=
\left( 0, d^{+}\alp + 2\pi f'_{\delta}(0)\omega_{0} 
\right)$$
where $\alp$ is some 1-form and $\omega_{0}:=\omega_{g_{0}}$ is as in
Proposition ~\ref{omegag}, and $g_{0}$ is the unique zero of the 
function $f_{\delta}$, as in ~\ref{derivative}.
\end{description}
\label{qderiv}
}
\end{lemma}
\vspace{.25in}
\noin
{\bf Proof:} Using the fact that 
$D_{A}~=~\sum c(e_{i})\circ\nabla_{A,e_{i}}$
and that $\nabla_{A+sa,e_{i}}=\nabla_{A,e_{i}} + 2\pi isa(e_{i})$, we
obtain that $D_{A+sa}=D_{A}+2\pi isa\circ (-)$. Thus :

$$\frac{d}{ds}_{|s=0}\left(D_{A+sa,t}(\Phi + s\phi)\right)= 
D_{A,t}\phi + 2\pi ia\circ\Phi $$
Finally, since $F_{A+sa}=F_{A}+s\,da$, we have $
F_{A+sa}^{+,t}=F_{A}^{+,t}+s\,d^{+}_{t}a$, so that 
$$\frac{d}{ds}_{|s=0}\left(\rho_{t}(F^{+,t}_{A+sa})\right)=\rho_{t}
(d_{t}^{+}a)$$
Skew-sesquilinearity of $\sigma$ implies that :
$$\frac{d}{ds}_{|s=0}\sigma(\Phi + s\phi, \Phi +s\phi)=
\sigma(\Phi,\phi)+\sigma(\phi,\Phi)=2\Im \sigma(\Phi,\phi)$$
and we have (i). 

To see (ii), since we are computing the derivative at $\Phi=0$,
in the direction of $\phi=0$, we have $D_{A,t}\Phi = 0$ and 
$\sigma (\Phi,\Phi)=0$ for all $t$, so the first coordinate
of the right hand side of the equation in the statement is clearly
zero. Now one just needs to compute $\frac{d}{dt}_{|t=0}
\left(\rho_{t}(F_{A}^{+,t})\right)$. By ~\ref{decomp}, we have 
$$\frac{d}{dt}_{|t=0}\left(\rho_{t}(F_{A}^{+,t})\right)
= \left\langle\frac{d}{dt}_{|t=0}
\left(\rho_{t}(F_{A}^{+,t})\right),\omega_{0}\right\rangle_{0}\omega_{0} 
+ d^{+}\alp$$
for some 1-form $\alp$, where $\langle\;,\;\rangle_{0}$ is the global inner product
on $\Omega^{2}$ with respect to $g_{0}$. Recall from \S~\ref{rhog} the 
facts that $\rho_{0}=\mbox{Id}$ and $(\rho_{t}(-),\rho_{t}(-))_{0}= 
(\;,\;)_{t}$, for the pointwise inner product, and so we have the same
formula for the global inner product, viz. 
$\langle\rho_{t}(-),\rho_{t}(-)\rangle_{0}= 
\langle\;,\;\rangle_{t}$. We compute the first 
term :
\begin{eqnarray*}
\langle\frac{d}{dt}_{|t=0}
\left(\rho_{t}(F_{A}^{+,t})\right),\omega_{0}\rangle_{0} &=&
\frac{d}{dt}_{|t=0}\langle\rho_{t}(F_{A}^{+,t}), 
\rho_{t}(\omega_{t})\rangle_{0}
-\langle F_{A}^{+,0},\frac{d}{dt}_{|t=0}\rho_{t}(\omega_{t})\rangle_{0}\\
= \frac{d}{dt}_{|t=0}\langle F_{A}^{+,t},\omega_{t}\rangle_{t}
-\frac{d}{dt}_{|t=0}\langle\delta, \rho_{t}(\omega_{t})\rangle_{0} &=&
\frac{d}{dt}_{|t=0}\left(\langle F_{A},\omega_{t}\rangle_{t} -
\langle\rho_{t}^{-1}\delta, \omega_{t}\rangle_{t}\right)\\
= \frac{d}{dt}_{|t=0}\left(2 \pi[c]\cup [\omega_t]- 
\langle(\rho_{t}^{-1}\delta)_{{\cal H}_{t}}, \omega_{t}\rangle_{t}\right)
&=& 2\pi \frac{d}{dt}_{|t=0}\left([c]\cup [\omega_{t}]-
\frac{1}{2\pi}[(\rho_{t}^{-1}(\delta))_{{\cal H}_{t}}]\cup 
[\omega_{t}]\right)\\
 = 2\pi \frac{d}{dt}_{|t=0}\left(\cdel\right)\cup[\omega_{t}]
&=& 2\pi f_{\delta}'(0)
\end{eqnarray*}
because for $\alp\in {\cal H}^{2}_{t}$, the $*_{t}$ self duality of
$\omega_{t}$ implies $[\alp]\cup [\omega_{t}]=\langle\alp,
 \omega_{t}\rangle_{t}$.
This proves the proposition.\hfill $\Box$

\subsection{Moduli Spaces}

\begin{proposition}{\rm 
There is a Baire subset $U$ of an $\eps$-ball around $0$ in $\Omega^{2+}$
for which both of the following hold : 
\begin{description}
\item[(i)] $(0,\delta)$ is a regular value for 
$Q_{|{\cal N}^{*}\times{\cal I}}$
\item[(ii)] The conclusions of Proposition~\ref{transverse1}, when
$b_{2}^{+}\geq 2$, and Proposition~\ref{derivative}, when $b_{2}^{+}=1$
are satisfied by $\delta$.
\end{description}
\label{transverse4}
}
\end{proposition}
\vspace{.25in}
\noin
{\bf Proof:} Since $Q$ is {\em not} a Fredholm map, one cannot directly
apply Proposition~\ref{fredmodel} or Corollary~\ref{transverse2} of the
Appendix to it. However, consider the space :
$$ {\cal M}^{*}({\cal I})
:=\{(A,\Phi,g_{t}):\;D_{A,t}\Phi=0,\;\Phi\neq 0,\;g_{t}\in 
{\cal I}\}$$
Since $Q$ is ${\cal G}$-equivariant, and ${\cal G}$ acts
trivially on ${\cal I}$ and $\Omega^{2+}$, it follows that the
${\cal G}$ action on ${\cal M}^{*}({\cal I})$ is free. Also, 
${\cal M}^{*}({\cal I})$ is fibred over ${\cal I}$ via
projection to the last coordinate. Note that 
${\cal M}^{*}({\cal I})= Q_{1}^{-1}(0)$ where 
$$Q_{1}=\pi_{\Gamma(W_{-})}\circ Q:{\cal N}^{*}\times {\cal I}\rarr
\Gamma(W_{-})$$
By \S 3.4 of [PP], we know that $0$ is a regular value of 
$Q_{1,t}:{\cal N}^{*}\times\{g_{t}\}~\rarr~\Gamma(W_{-})$. Thus
it is a regular value of $Q_{1}$, so that 
${\cal M}^{*}({\cal I})$ is an (infinite dimensional)
submanifold of ${\cal N}^{*}\times{\cal I}$, whose fibre over 
$g_{t}\in {\cal I}$ is ${\cal M}^{*}(g_{t})$ (as defined in
\S~3.4 of [PP]). 
Now, by Lemma~\ref{transverse3} of the Appendix, 
regular values $(0,\delta)$
will exist for $Q_{|{\cal N}^{*}\times{\cal I}}$ whenever
regular values exist for the map :
\begin{eqnarray*}
Q_{2}=\pi_{\Omega^{2+}}\circ Q:{\cal M}^{*}({\cal I})&\rarr &\Omega^{2+}\\
(A,\Phi,g_{t})&\mapsto &\left(\rho_{t}(F_{A}^{+,t})-\sigma (\Phi,\Phi)\right)
\end{eqnarray*}
Again, $Q_{2}$ is {\em not} a Fredholm map. On the other hand, it is constant
along ${\cal G}$-orbits. 

\vspace{.25in}
\noin
{\bf Claim 1:} $Q_{2}$ descends to a map 
${\cal M}^{*}({\cal I})/{\cal G}~\rarr~\Omega^{2+}$ (which we
also denote by $Q_{2}$), and this map is Fredholm of index $d(L)+1$
where 
$$d(L):=\frac{1}{4}\left(c_{1}(L)^{2}-3\sigma(X)-2\chi(X)\right)$$
($\sigma(X)$ is the signature, and $\chi(X)$ the Euler characteristic
of $X$).

\vspace{.25in}
\noin
{\bf Proof of Claim 1:} Since $\rho_{t}(F_{A}^{+,t})~-
~\sigma (\Phi, \Phi)=
\rho_{t}(F_{gA}^{+,t})~-~\sigma (g\Phi, g\Phi)$ 
for all $g\in {\cal G}$, it is clear that $Q_{2}$ descends to
the quotient space ${\cal M}^{*}({\cal I})/{\cal G}$. 

Now, the inclusion of tangent spaces: $T_{x}({\cal M}^{*}(g_{t})/{\cal G})
\hookrightarrow
T_{x}({\cal M}^{*}({\cal I})/{\cal G})$ has codimension one at
each $x$ (= dimension of the tangent space to ${\cal I}$), it 
is enough to prove that 
$$Q_{2,t}:{\cal M}^{*}(g_{t})/{\cal G}~\rarr~\Omega^{2+}$$ 
is a Fredholm map of index $d(L)$. The diagram :
$$\begin{array}{cclclclcc}
0&\rarr&T_{e}({\cal G})&\longrightarrow &T_{(A,\Phi,g_{t})}
({\cal M}^{*}(g_{t}))&\longrightarrow
&T_{([A,\Phi ],g_{t})}({\cal M}^{*}(g_{t})/{\cal G})&\rarr &0\\
& &\downarrow\,d^{*}d=\Delta& &\downarrow\,\chi& &\downarrow\,DQ_{2,t}& & \\
0&\rarr
&\Omega^{0}&\hookrightarrow&\Omega^{2+}\oplus\Omega^{0}
&\longrightarrow &\Omega^{2+}&\longrightarrow & 0
\end{array}
$$
shows that $DQ_{2,t}$ on the right is Fredholm of index $d(L)$
iff the middle map $\chi$ given by $\chi (a,\phi, 0)=(DQ_{2,t}(a,\phi),d^{*}a)$ 
is Fredholm of index $d(L)$,
since the laplacian $\Delta$ has one-dimensional kernel and
cokernel, $X$ being connected. 

However, the diagram :

$$\begin{array}{cclclclcc}
0&\rarr& T_{(A,\Phi,g_{t})}({\cal M}^{*}(g_{t})) &\rarr &T_{(A,\Phi,g_{t})}
({\cal N}^{*}\times \{g_{t}\})&\stackrel{DQ_{1,t}}{\rarr} 
&\Gamma(W_{-})&\rarr &0\\
& &\downarrow\,\chi& &\downarrow\,\chi_{1}=(DQ_{t},d^{*})& &\left|\!\right|& & 
\\
0&\rarr
&\Omega^{2+}\oplus\Omega^{0}&\hookrightarrow&\Gamma (W_{-})\oplus\Omega^{2+}
\oplus\Omega^{0}&\longrightarrow &\Gamma(W_{-})&\rarr & 0
\end{array}
$$
where $\chi_{1}(a,\phi, 0)=(DQ_{t}(a,\phi), d^{*}a)$, 
implies that $\chi$ is a Fredholm operator iff $\chi_{1}$ is a
Fredholm operator, and of the same index since the vertical map
on the right is an equality. Now, note that by the Proposition~\ref{qderiv},
part (i), we have :
$$\chi_{1}(a,\phi,0)=(D_{A,t}\phi + 2\pi ia\circ\Phi, \rho_{t}(d_{t}^{+}a)-
2\Im(\sigma (\Phi,\phi)),d^{*}a)$$
Now, $a\in L^{2}_{k},\;\Phi\in L^{2}_{k}$ implies that
$a\circ\Phi\in L^{2}_{k}$, (Leibnitz rule for k-th derivative of
a product and Schwartz inequality), and since we have
$L^{2}_{k-1}$ Sobolev norm on $\Gamma(W_{-})$, Rellich's Lemma implies
that the map $a~\mapsto~a\circ\Phi$ is a compact operator from 
$\Omega^{1}_{L^{2}_{k}}$ into $\Gamma(W_{-})$. Similar considerations
apply to $\Im\sigma(\Phi,\phi)$. Thus $\chi_{1}$ is a compact 
perturbation of the map :
\begin{eqnarray*}
\Gamma(W_{+})\times \Omega^{1}&\rarr &\Gamma(W_{-})\times \Omega^{2+}
\times\Omega^{0}\\
(a,\phi)&\mapsto & (D_{A,t}\phi, \rho_{t}(d_{t}^{+}a), d^{*}a)
\end{eqnarray*}
whose index is clearly $\mbox{index}(D_{A})~+~\mbox{index}(d^{+},d^{*})$.
The index of the second map is the negative of the 
Euler characteristic of the complex
$$\Omega^{0}\stackrel{d}{\rarr}\Omega^{1}\stackrel{d^{+}}{\rarr}
\Omega^{2+}$$
which is $\left(-\dim H^{2+}-\dim H^{0} + \dim H^{1}\right)$, which is
$-\frac{1}{2}\left(\sigma(X)+\chi(X)\right)$. The Atiyah-Singer
index theorem for the Dirac operator implies :
$$\mbox{index}D_{A} = -\frac{1}{4}\left({\sigma(X)} - c_{1}(L)^{2}\right)
$$
Combining these, we have the index of 
$Q_{2}~:~{\cal M}^{*}({\cal I})/{\cal G}~\rarr~\Omega^{2+}$ to be 
$$d(L):=\frac{1}{4}\left(c_{1}(L)^{2}-3\sigma(X)-2\chi(X)\right)$$
which proves Claim 1. \hfill $\Box $

We now return to the proof of our Proposition. By ~\ref{transverse2}
and ~\ref{transverse3}, we have a Baire subset of a
neighbourhood of $0$ in $\Omega^{2+}$ such that for $\delta$ in this
subset, $(0,\delta)$ is a regular value for 
$Q_{|{\cal N}^{*}\times{\cal I}}$. This proves (i). To get (ii),
one just intersects this Baire subset with the $\eps$-ball $U$ that
is guaranteed by Propositions~\ref{transverse1} and ~\ref{derivative}.
\hfill $\Box$

\begin{notation}{\rm
We now fix a $\delta$ so that both (i) and (ii) of the last 
 Proposition~\ref{transverse4} are satisfied. In the case $b_{2}^{+}=1$,
we assume that the function $f_{\delta}:{\cal I}\rarr \R$ has its unique
zero at $g_{0}$, and $0$ is a regular value for it.
 For notational simplicity we modify $Q$ to 
$Q_{\delta}$, a translate of $Q$, by the formula :
$$Q_{\delta}(A,\Phi, g_{t})=Q(A,\Phi, g_{t}) - (0,\,\delta)$$
We thus have the following consequence to the propositions and 
corollaries ~\ref{fredmodel}, ~\ref{transverse2},
~\ref{transverse3},~\ref{transverse3} and ~\ref{transverse4}:  }
\label{notation1}
\end{notation}

\vspace{.25in}
\begin{proposition}{\rm We have the following facts about $Q_{\delta}$:
\begin{description}
\item[(i)] The map :
$$Q_{\delta}:{\cal N}\times {\cal I}\rarr \Gamma(W_{-})\times\Omega^{2+}$$
is a ${\cal G}$-equivariant map with $(0,0)$ as a regular value for
$Q_{\delta |{\cal N}^{*}\times {\cal I}}$, so that 
$Q_{\delta}^{-1}(0,\,0)\cap {\cal N}^{*}\times {\cal I}$ is a
Banach manifold, which we will denote as ${\cal
M}_{\delta}^{*}({\cal I})$. It is fibred over ${\cal I}$ with
fibre ${\cal M}_{\delta}^{*}(g_{t})$ on $g_{t}$.

\item[(ii)] When $b_{2}^{+}\geq 2$, $Q_{\delta}^{-1}(0,\,0)\subset
{\cal N}^{*}\times{\cal I}$. In this case the group of gauge
transformations ${\cal G}$ acts {\em freely} on all of 
$Q_{\delta}^{-1}(0,\,0)$, and consequently the quotient space
$Q_{\delta}^{-1}(0,\,0)/{\cal G}={\cal M}_{\delta}^{*}({\cal
I})/{\cal G}:=M_{c,\delta}({\cal I})$ is a  manifold of 
dimension $d(L)+1$. This manifold may also be regarded as
  $Q_{2,\delta}^{-1}(0)$, where
$Q_{2,\delta}=\pi_{\Omega^{2+}}\circ Q_{\delta}$ is as in Claim 1
in the proof of Proposition~\ref{transverse4}. Its fibre over $g_{t}$,
 $M_{c,\delta}(g_{t})$ is a {\em compact } manifold of dimension $d(L)$
for all $g_{t}\in {\cal I}$. Thus $M_{c,\delta}({\cal I})$ is
also compact. In particular, when $d(L)=0$, it is a finite union of 
arcs, and the cardinality $\#(M_{c,\delta}(g_{t}))\;(\mbox{mod}\;2)$
is independent of $t$, and consequently an invariant of $X$
(with its given $\mbox{Spin}_{c}$ structure). 
\item[(iii)] If $b_{2}^{+}=1$, and $H^{1}(X,\R)=0$, the space

$$\left(Q_{\delta}^{-1}(0,\,0)\cap ({\cal A}\times\{0\}\times
{\cal I})\right)/{\cal G}$$ 
is just a single point (called a {\em reducible solution}), say
 $([A_{0},0],g_{0})$. A 
neighbourhood of this point in $Q_{\delta}^{-1}(0,\,0)/{\cal G}$
is homeomorphic to $\phi^{-1}(0,\,0)/S^{1}$ where :
$$\phi : \kernel{D_{A_{0}, 0}}\rarr \coker{D_{A_{0},0}} $$
is a smooth map which is (a) $S^{1}$-equivariant, and (b) has $0$
as a regular value in a small deleted neighbourhood $U -\{0\}$ of
$0$ in $\kernel{D_{A_{0},0}}$. For metrics $g_{t}$ such that $t\neq 0$
the moduli space $M_{c,\delta}(g_{t})$ is a finite set of points if
$d(L)=0$.

\end{description}
}
\label{moduli1}
\end{proposition}
\vspace{.25in}
\noin
{\bf Proof:} 

(i) follows immediately from Proposition~\ref{transverse4}.
For (ii), note that $Q_{\delta}(A,\Phi,g_{t})=(0,\,0)$ implies that
$D_{A,t}\Phi~=~0$, and
$\rho_{t}(F_{A}^{+,t})=\sigma(\Phi,\Phi)+\delta$.
This implies 
 
$$\frac{1}{2\pi}\left(\rho_{t}^{-1}\sigma(\Phi,\Phi)\right)_{{\cal H}_{t}}=
\frac{1}{2\pi}\left(F_{A}^{+,t}-
\rho_{t}^{-1}(\delta)\right)_{{\cal H}_{t}}=
\pi^{+}_{t}\circ\psi_{t}(c(g_{t},\delta))=
\sigma_{c,\delta}(g_{t})$$
in the notation of ~\ref{transnotation}.
But ${\cal I}\cap B_{c,\delta}=\phi$, i.e. the section $\sigma_{c,\delta}$
is nonvanishing on ${\cal I}$ by our choice of ${\cal I}$ and
$\delta$, so $\sigma(\Phi,\Phi)\neq 0$, and so $\Phi\neq
0$. Thus there are no reducible solutions, and 
$Q_{\delta}^{-1}(0,0)\subset {\cal N}^{*}\times{\cal I}$.  Hence 
$Q_{\delta}^{-1}(0,0)/{\cal G}$ is a smooth manifold, since the ${\cal
G}$ action is free on ${\cal N}^{*}\times {\cal I}$, and since each
fibre $M_{c,\delta}(g_{t})$ (notation of \S 3.5 of [PP]) is compact by
\S 5.2 of [PP]. Since ${\cal I}$ is compact, so is $M_{c,\delta}({\cal
I})$.  In the case when $d(L)=0$,
$M_{c,\delta}(g_{t})$ is then a finite set of points, and
$M_{c,\delta}(g_{a})$ and $M_{c,\delta}(g_{b})$ are {\em cobordant},
and so have the same cardinality (modulo 2). In particular,
$M_{c,\delta}(g_{1})$ and $M_{c,\delta}(g_{-1})$ have the same
cardinality modulo 2. Thus we have (ii), since the rest of it follows
from Proposition~\ref{transverse4} and Lemma~\ref{transverse3} of the
Appendix.

To get (iii), note that when 
$b_{2}^{+}=1$, $Q_{2,\delta}(A,0,g_{t})=(0,\,0)$ implies
that $\sigma_{c,\delta}(g_{t})=0$. This happens (see 
Proposition~\ref{derivative}) iff $f_{\delta}(g_{t})=0$. By 
 our choice of ${\cal I}$ and $\delta$, this happens only at $g_{t}=g_{0}$. 
Now choose
a point $(A_{0},0,g_{0})\in Q_{\delta}^{-1}(0,\,0)\cap ({\cal A}\times
\{0\}\times{\cal I})$. Then $Q_{\delta}(A_{0}+a,0,g_{0})=(0,\,0)$
implies $\rho_{0}(d_{0}^{+}a)=d^{+}a=0$. Let $g\in{\cal G}$. The gauge
action takes $A_{0}$ to $A_{0}+g^{*}\omega$, where $\omega\in 
H^{1}(S^{1},\R)$ is the generating ``angle'' 1-form of $S^{1}$. But
$g^{*}\omega$ is thus a closed 1-form on $X$. Since
$H^{1}(X,\R)=0$ by assumption, $g^{*}\omega=d\alp$ for some function
$\alp\in\Omega^{0}$. Conversely, given an $\alp\in \Omega^{0}$,
the map $g:X\rarr S^{1}$ defined by $g(x)=e^{2\pi i\alp (x)}$ satisfies
$g^{*}(\omega)=d\alp$. Thus 
$\left(Q_{\delta}^{-1}(0,\,0)\cap ({\cal A}\times
\{0\}\times{\cal I})\right)/{\cal G} \simeq \frac{\kernel{d^{+}}}{\im{d}}
=H^{1}(X,\R)=0$. So it consists of a single
point $[(A_{0}, 0, g_{0})]$.

Now we need to get a model for a neighbourhood of this point. It is well 
known that for a smooth action, the neighborhood of a point in the orbit
space corresponding to an orbit with isotropy $G$ is homeomorphic to a
neigborhood in the orthogonal slice of that orbit divided by the isotropy 
$G$. A slice in $\Omega^{1}\times\Gamma(W_{+})\times{\cal I}$
orthogonal to the ${\cal G}$-orbit of $(A_{0}, 0,g_{0})$, which we
will take as $(0,0,g_{0})$ by setting the origin at $A_{0}$, is
clearly $(\im{d})^{\perp}\times\Gamma(W_{+})\times\R$ because the 
${\cal G}$-orbit has tangent space $\im{d}\times 0\times 0$ at 
$(0,0,g_{0})$. Also $DQ_{\delta\,;A_{0},0,g_{0}}$ is an
$S^{1}$-equivariant map, (because $Q_{\delta}$ is ${\cal G}$
equivariant), where $S^{1}$ is the isotropy group of
the point $(A_{0},0,g_{0})$, and of course, this $S^{1}$-action
is orthogonal and linear (it is the derivative of the
$S^{1}$-action on the space 
${\cal A}\times\Gamma(W_{+})\times{\cal C}$). By
Lemma~\ref{qderiv}, we have :
\begin{eqnarray*}
DQ_{\delta\,;A_{0},0,g_{0}}(a,\phi,0)&=&(D_{A_{0}}\phi , d^{+}a)\\
DQ_{\delta\,;A_{0},0,g_{0}}(0,0,\lambda g'(0))
&=&(0,\lambda(2\pi f_{\delta}'(0)
\omega_{g_{0}}+ d^{+}\alp))
\end{eqnarray*}
We now claim that:
$$DQ_{\delta\,;A_{0},0,g_{0}}:
(\im{d})^{\perp}\times\Gamma(W_{+})\times\R
\rarr\Gamma(W_{-})\times\Omega^{2+}$$
is a Fredholm operator, whose kernel is 
$\kernel{D_{A_{0}}}\subset\Gamma(W_{+})$, and 
cokernel is $\coker{D_{A_{0}}}~\subset~\Gamma(W_{-})$. 

From the above formulae it follows that
$$\kernel{(DQ_{\delta\,;A_{0},0,g_{0}})}=
\left\{(a,\phi,\lam g'(0)): a\,\perp\im{d},\,D_{A_{0}}\phi = 0,
d^{+}a +\lam(d^{+}\alp + 2\pi f_{\delta}'(0)\omega_{g_{0}})=0\right\}$$
Now, $\omega_{g_{0}}\perp \;\im{d^{+}}$, so on the right hand
side we must have $\lam 2\pi f'_{\delta}(0)\omega_{g_{0}}=0$. But
$f_{\delta}'(0)\neq 0$, by our choice of ${\cal I}$, $\delta$ and 
Proposition~\ref{derivative},
so $\lam =0$. But this implies that $D_{A_{0}}\phi=0$ and $d^{+}a=0$.
Thus 
$$\kernel{(DQ_{\delta\,;A_{0},0,g_{0}})}=\left((\im{d})^{\perp}\cap
\kernel{d^{+}}\right)\times \kernel{D_{A_{0}}}$$
Since $(\im{d})^{\perp}\cap\kernel{d^{+}}\simeq \kernel{d^{+}}/\im{d}
\subset H^{1}(X,\R)=0$, this is just $\kernel{D_{A_{0}}}$. It is
finite dimensional by ellipticity of the Dirac operator.

To show that 
$$DQ_{\delta\,;A_{0},0,g_{0}}((\im{d})^{\perp}\times\Gamma(W_{+})
\times \R)$$
is closed, it is enough to show that 
$DQ_{\delta\,;A_{0},0,g_{0}}((\im{d})^{\perp}\times\Gamma(W_{+})
\times 0)$ is closed. But this is just $\im{D_{A_{0}}}\times d^{+}
(\im{d})^{\perp}$. Now $(\im{d})^{\perp}=\kernel{d^{*}}={\cal H}^{1}
\oplus d^{*}\Omega^{2}$, so $d^{+}(\im{d})^{\perp}=d^{+}d^{*}\Omega^{2}$.
Since $\Omega^{1}={\cal H}^{1}\oplus d^{*}\Omega^{2}
\oplus d\Omega^{0}$, $d^{+}\Omega^{1}=d^{+}d^{*}\Omega^{2}$ as
well. Thus $d^{+}(\im{d})^{\perp}=d^{+}\Omega^{1}$. This is
clearly closed by the decomposition of Lemma~\ref{decomp}. Since
$D_{A_{0}}$, the Dirac operator, is elliptic, its range is
closed too, so the range of 
$$DQ_{\delta\,;A_{0},0,g_{0}}:
(\im{d})^{\perp}\times\Gamma(W_{+})\times\R\rarr\Gamma(W_{-})\times\Omega^{2+}$$
is closed. Its cokernel is :
$$\{(\psi,\,\tau):\psi\perp\,\im{D_{A_{0}}};\;\tau\perp 
d^{+}a +\lam(d^{+}\alp + 2\pi f_{\delta}'(0)\omega_{0}\;\forall\,\lam\in\R,
\,a\in(\im{d})^{\perp}\}$$
This implies $\psi\in\coker{D_{A_{0}}}$. Also, setting $\lam =0$,
one finds that $\tau\perp d^{+}((\im{d})^{\perp})$, i.e. $\tau\perp
d^{+}\Omega^{1}$. Since $\tau\in\Omega^{2+}$, the Lemma~\ref{decomp}
implies that $\tau\in{\cal H}^{2+}$. However, by setting $a=0$ and
$\lam =1$, $\tau$ is also orthogonal to 
$2\pi f_{\delta}'(0)\omega_{0}+
d^{+}\alp$, and since we already have $\tau \perp d^{+}\Omega^{1}$,
it follows that $\tau\perp\omega_{0}$, since $f_{\delta}'(0)\neq 0$
by our choice of ${\cal I}$ (from~\ref{derivative}). We now note
that $\omega_{0}$ was chosen as the basis element of ${\cal H}^{2+}$,
so $\tau=0$. Thus the cokernel of our map is just $\coker{D_{A_{0}}}$.
This proves our Fredholm-ness assertion, and the identifications of
kernel and cokernel.

Thus, Proposition~\ref{equivmodel2} applies, and a neighbourhood of 
$(A,0,g_{0})$
in $Q_{\delta}^{-1}(0,0)/{\cal G}$ is homeomorphic to a
neighbourhood of $0$ in $\phi^{-1}(0)/S^{1}$, where the $S^{1}$-equivariant
 map :

$$\phi :\kernel{D_{A_{0}}}\rarr \coker{D_{A_{0}}}$$
has $0$ as a regular value when restricted to
$\phi^{-1}(0)-\{0\}$. 
The last statement follows from the fact that if $t\neq 0$ then
$f_{\delta}(g_{t})\neq 0 $,
and ${\cal M}^{*}_{\delta}(g_{t})={\cal M}_{\delta}(g_{t})$ is a
Banach manifold (by part (i)), since there are no reducible solutions,
and its quotient by ${\cal G}$ has dimension $d(L)=0$ by part (ii) above. 
This proves (iii), and the proposition.\hfill $\Box$

\section{Computations}
\subsection{The case of $\C\Proj^{2}\# n\overline{\C\Proj^{2}}$}
\label{peetwo}
Let $X=\C\Proj^{2}\# n\overline{\C\Proj^{2}}$, which is just the
blow-up of $\C\Proj^{2}$ at $n$ points. The case $n=0$ is just $\C\Proj^{2}$,
which is a special case of the ensuing discussion. Since $X$ is a 
complex manifold, we have a canonical $\mbox{spin}_{c}$ structure on $X$ 
(as in \S~6.1 of [PP]), with $L=K_{X}^{-1},\;W_{+}=K_{X}^{-1}\oplus
{\bf 1}_{\C},\;W_{-}=T_{hol}(X)$. Also 
\begin{description}
\item[(i)] $H^{1}(X,\R)=0$, and
\item[(ii)] $H^{2}(X,\R)=\R H \oplus \sum_{i=1}^{n}\R E_{i}$,
where $H$ is the pullback of the hyperplane class in
$\C\Proj^{2}$ via the blow up map $\pi$ and will, from here on, be 
called the hyperplane class by abuse of language. $E_{i}$ is the
generating $\overline{\C\Proj}^{1}$ in the $i$-th copy of 
$\overline{\C\Proj}^{2}$. The cup products between these classes
are :
$$H\cup H=1,\;E_{i}\cup E_{j}=-\delta_{ij},\;H\cup E_{i}=0\;\;\forall\;\;1\leq
i, j\leq n$$
Thus $b_{2}^{+}=1,\;b_{2}^{-}=n$ and the cup pairing is of type
$(1,n)$. 
\item[(iii)] From the formula
$K_{X}=\pi^{*}(K_{\C\Proj^{2}})\otimes [E]$, where
$E:=\sum_{i=1}^{n} E_{i}$ is the exceptional divisor, and that 
$c_{1}(K^{-1}_{\C\Proj^{2}})=c_{1}(\C\Proj^{2})=3H_{\C\Proj^{2}}$, it
follows that $c_{1}(L)=c_{1}(K_{X}^{-1})=3H - E$. Thus $c_{1}(L)^{2}=9-n$,
and this is $< 0$ whenever $n > 9$. Since $\chi(X)=n+3$, and 
$\sigma(X)=1-n$, we have
$d(L)=\frac{1}{4}\left(c_{1}(L)^{2}-3\sigma -2\chi \right)=0$ (see 
Proposition~\ref{transverse4} for the definition of $d(L)$).
\end{description}

\vspace{.25in}
\noin
\begin{proposition}{\rm
For $X=\C\Proj^{2}\# n\overline{\C\Proj^{2}}$,
the Seiberg-Witten Moduli Space $M_{c,\delta}(g_{t})$ consists
of finitely many points whenever $f_{\delta}(g_{t})\neq 0$, 
(i.e. when $t\neq 0$) where $f_{\delta}$ is the
function defined in ~\ref{derivative}. 
}
\label{blowup}
\end{proposition}

\vspace{.25in}
\noin
{\bf Proof:} Follows immediately from (iii) of Proposition~\ref{moduli1}
\hfill $\Box$

\vspace{.25in}
Of course, this proposition does not tell us what the cardinality
(mod 2) of the moduli space $M_{c,\delta}(g_{t})$ might be. To show
that there exist metrics for which this cardinality is non-zero
will be our next goal. 

\begin{proposition}{\rm Let $X=\C\Proj^{2}\#
n\overline{\C\Proj^{2}}$, as above, with $n > 9$.
Let ${\cal I}$ be an arc in ${\cal C}$ chosen in accordance with
the Proposition~\ref{moduli1} (i.e. $(0,0)$ is a regular value
of $Q_{\delta |{\cal N}^{*}\times {\cal I}}$, and 
the function $f_{\delta}:{\cal I}\rarr\R$  satisfies
the conclusion of Proposition~\ref{derivative}. Then :
$$\#M_{c,\del}(g_{-1}) -\#M_{c,\del}(g_{1})=1\;\;(\mbox{mod}\,2)$$
}
\label{moduli2}
\end{proposition}
\vspace{.25in}
\noin
{\bf Proof:} By (iii) of Proposition~\ref{moduli1}, a neighbourhood of $(A_{0},0,g_{0})$,
the unique reducible point in $M_{c,\delta}({\cal I})$ is
homeomorphic to a neighbourhood of $0$ in $\phi^{0}/S^{1}$, where
$\phi :\kernel{D_{A_{0}}}\rarr\coker{D_{A_{0}}}$ is a smooth 
$S^{1}$-equivariant map, and $\phi$ has $0$ as a regular value in
a deleted neighbourhood of $0$ in $\kernel{D_{A_{0}}}$. We need
to show that an odd number of arcs emerge from $0$ in this neighbourhood.

Now
\begin{eqnarray*}
\mbox{index}_{\C}(D_{A_{0}})&=&\frac{1}{2}\mbox{index}_{\R}(D_{A_{0}})
= \frac{1}{8}(c_{1}(L)^{2}-\sigma(X))\\
&=&\frac{1}{8}(9-n-1+n)=1
\end{eqnarray*}
Thus if $\dim_{\C}(\coker{D_{A_{0}}})=r$, then 
$\dim_{\C}(\ker{D_{A_{0}}})=r+1$.
So our $\phi$ is an $S^{1}$-equivariant smooth map (with $S^{1}$
acting as scalar multiplication on both sides) from $\C^{r+1}$
to $\C^{r}$, with $0$ a regular value for $\phi_{|\C^{r+1}-0}$. Let
${\cal O}(-1)$ denote the tautological bundle on $\C\Proj^{r}$,
and let 
$$\pi:(\C^{r+1}-0)/S^{1}\simeq\C\Proj^{r}\times\R_{+}\rarr \C\Proj^{r}$$
denote the projection into the first factor. If we then denote
by $<z_{0},...z_{r}>$ the $S^{1}$-equivalence class of $(z_{0},..z_{r})$
in the orbit space $(\C^{r+1}-0)/S^{1}$, we get a natural section
$s$ of
the bundle $\pi^{*}\hom ({\cal O}(-1),\C^{r})\simeq\pi^{*}
(\C^{r}\otimes {\cal O}(1))$ by setting :
$$s(<z_{0},...,z_{r}>)(u)=\phi(z_{0},...,z_{r})$$
where $u$ is the unit vector $\frac{(z_{0},...,z_{r})}{\norml(z_{0},..,z_{r})
\normr}$. The $S^{1}$-equivariance of $\phi$ implies that this
map $s$ is well defined, and that $0$ is a regular value for $\phi$
on a deleted neighbourhood means that $s$ is transverse to the zero-section. 
Now, $s^{-1}(0)$ is clearly $\phi^{-1}(0,0)-\{0\}/S^{1}$. Thus, modulo 2,
the number of arcs going to $0$ is precisely the number of
points in $s^{-1}(0)\cap (\C\Proj^{r}\times\{\eps\})$ modulo 2, for generic
$\epsilon$. But this is just the Euler number of the bundle 
$\C^{r}\otimes{\cal O}(1)$, which is 1. Thus 
$$\#M_{c,\del}(g_{-1}) -\#M_{c,\del}(g_{1})=1\;\;(\mbox{mod}\,2)$$
and we are done.\hfill $\Box$

\vspace{.25in}
\noin
\begin{proposition}{\rm(Hitchin) There exist metrics $g$ on 
$X=\C\Proj^{2}\# n\overline{\C\Proj^{2}}$, such that :
\begin{description}
\item[(i)] $g$ is Kahler.
\item[(ii)] $[\omega_{g}]\cup c_{1}(L)$ has the same sign as the scalar
curvature $s_{g}$, which is positive.
\item[(iii)] For any metric $g'$ conformally equivalent to a Kahler metric
$g$ satisfying (i) and (ii) above, $[\omega_{g'}]\cup c_{1}(L) > 0$. 
\end{description}
}
\label{hitchin}
\end{proposition}

\vspace{.25in}
\noin
{\bf Proof:} See reference number [4] in [KM] for the proof of (i) and
(ii). For (iii), note that $[\omega_{g'}]=[\omega_{g}]$ when $g$ and $g'$
are conformally equivalent, by (iii) of ~\ref{omegag}. \hfill$\Box$

\vspace{.25in}
Clearly, for a metric as in Proposition~\ref{hitchin} above, the
moduli space $M_{c,\delta}(g)$ is empty, by \S~5.2 of [PP]. 
Consequently, we have the following corollary to Proposition~\ref{moduli2}.

\vspace{.25in}
\begin{corollary}{\rm If $g$ is a metric on 
$X~=~\C\Proj^{2}\# n\overline{\C\Proj^{2}}$ such that $c_{1}(L)
\cup [\omega_{g}]~<~0$, then the moduli space $M_{c,\delta}(g)\neq\phi$.
}
\label{nonempty}
\end{corollary}

\vspace{.25in}
\noin
{\bf Proof:} Assuming there is such a metric, join it by an arc
${\cal I}$ in ${\cal C}$ meeting $B_{c}$ transversely, at one point, 
and apply Propositions ~\ref{derivative}, \ref{moduli2}.\hfill $\Box$

\subsection{The Tubing Construction}
Let $X$ be a compact, connected oriented Riemannian 4-manifold,
and let $Y$ be a compact 3-manifold, also oriented (so that its normal
bundle in $X$ is trivial) such that :
$$X - Y = X_{+} \cup X_{-}$$
as two disjoint components. Assume further that there exists a
Riemannian metric $g$ on $X$ such that 
$g_{|\nu (\eps)}=dt^{2}\times g_{Y}$, where $\nu (\eps)$ is an $\eps$-
tubular neighbourhood of $Y$ in $X$, and $g_{Y}$ is a smooth
Riemannian metric on $Y$. Let us denote $g_{|X_{\pm}}:=g_{\pm}$.
\begin{definition}{\rm Define the {\em R-tubing of $X$}, denoted
$(X_{R}, g_{R})$ to be the manifold.

$$X_{R}=(X_{-}\cup ([-R,0]\times Y))\cup (([0,R]\times Y)\cup X_{+})$$
with the metric $g_{R}$ being defined by $g_{R|X_{\pm}}=g_{\pm}$
and $g_{|[-R,R]\times Y}= dt^{2}\times g_{Y}$ on 
the piece $[-R,R]\times Y$. Note the ends $\{\pm R\}\times Y$ are
identified with $\partial X_{\pm}$. 
}    
\label{tubing}
\end{definition}
Clearly, since $X_{R}$ is diffeomorphic to $X$, one may view $g_{R}$ as 
a new metric on $X$. Of course, its volume form $dV_{g_{R}}$ will change,
but by a global conformal change one may restore the old volume form. 
Note that this conformal change does not affect $\omega_{0}$. Now let
us go back to $X=\C\Proj^{2} \# d^{2}\overline{\C\Proj}^{2}$, where 
$d > 3$ (so that $n=d^{2} > 9$ and $c_{1}(L)^{2}=9 - n < 0$ from the
opening discussion of this section). 

Let $\Sigma $ be an oriented smoothly embedded surface in
$\C\Proj^{2}$ whose homology class is Poincare dual to $dH$ (i.e., it
is of degree $d$) in $H^{2}(\C\Proj^{2})$. Thus the cap (or
Kronecker) product $H.[\Sigma]=d$. 

Let $S_{i}$ be a sphere $\C\Proj^{1}$ in the $i$-th copy of 
$\overline{\C\Proj}^{2}$ which is dual to $(-E_{i})$, so that 
$E_{j}.S_{i}~=~-\del_{ij}$. Let us  define the internal connected sum 
$\til{\Sigma}=\Sigma \# S_{1}\# S_{2}..\# S_{d^{2}}$. Clearly the genus 
$g(\til{\Sigma})$ is the same as that of $\Sigma$, i.e. $g(\Sigma)$, since 
$\Sigma$ and $\til{\Sigma}$ are homeomorphic. Also $[\til{\Sigma}]$ is
the homology class Poincare dual to $dH - E$. 

\begin{lemma}{\rm 
$\til{\Sigma}$ has trivial normal bundle in $X$.
}
\label{trivialnormal}
\end{lemma}

\vspace{.25in}
\noin
{\bf Proof:} Let $\nu =\nu(\til{\Sigma})$ be the normal bundle of
$\til{\Sigma}$ in $X$. Then we know that the Euler number of this
bundle is precisely the self-intersection number of $\til{\Sigma}$ in
$X$, i.e. $(dH -E).(dH - E)=d^{2}H^{2}+E^{2}=d^{2}+\sum
E_{i}^{2}=d^{2}-d^{2}=0$. Thus $\nu(\til{\Sigma})$ has a nowhere
vanishing section, i.e. has a trivial line bundle as a summand. Since
it is orientable, it follows that it is trivial. \hfill $\Box$

\vspace{.25in}

Thus if $\nu_{\eps}(\til{\Sigma})$ 
is an $\eps$-tubular neighbourhood
of $\til{\Sigma}$, we have 
$\nu_{\eps}(\til{\Sigma})\simeq D_{\eps}\times\til{\Sigma}$, and its
boundary $\partial\nu_{\eps}(\til{\Sigma})\simeq
S^{1}\times\til{\Sigma}$. Call this last space $Y$. Now apply the
tubing construction of ~\ref{tubing} to $X, \,Y$, with $X_{-}= 
X -\nu_{\eps}(\til{\Sigma}),\,X_{+}=\nu_{\eps}(\til{\Sigma}),\,
Y=\partial\nu_{\eps}(\til{\Sigma})$. Note that an $\eps$-neighborhood
of $Y$ is diffeomorphic to $(-\eps,\eps)\times Y$. 

\begin{proposition}{\rm
Assume there exists a metric $g_{0}$ on $X$ such that $g_{0|\nu_{\eps}(Y)}$ is
a product metric, and let $(X(R),g_{R})$ denote the $R$-tubing
of $g$ as defined above. Then,
for $R$ sufficiently large, $c_{1}(L)\cup [\omega_{g_{R}}] < 0$. 
}
\label{negativecup}
\end{proposition}

\vspace{.25in}
\noin
{\bf Proof:} For notational ease we shall denote $\omega_{g_{R}}$ by
$\omega_{R}$. Let $R_{i}$ be a sequence of positive numbers tending to
$\infty$. We know by Proposition~\ref{omegag} that in
$H^{2}(X,\R)={\cal H}$,   $[\omega_{R_{i}}]\cup [\omega_{R_{i}}]=1 $, and 
$[\omega_{R_{i}}]\in C_{+}$ where
$C_{+}$ is the preferred component of the positive cone of the
cup-product form, which contains the class $H$ defined above. So
$[\omega_{R_{i}}]\cup H > 0$ for all $i$. Let us normalise and define
$\omega_{i}$ to be the unique $\Delta_{0}$ harmonic 2-form
representing $( [\omega_{R_{i}}]\cup H )^{-1}[\omega_{R_{i}}]$. It is
clearly enough to show that $c_{1}(L)\cup\omega_{i}~<~0$ for $i$ 
large enough. 
Since
$\omega_{i}\cup H = 1$ for all $i$, and hence $\omega_{i}$ lie in the
bounded region $C_{+}\cap \{\bet:\,\bet\cup H = 1\}$, so have a
 (with respect to the global inner product
$\langle\,,\,\rangle_{g_{0}}$) $L^{2}$-convergent subsequence, which we continue to call
$\omega_{i}$. Since $\omega_{i}$ are $\Delta_{g_{0}}$-harmonic, the
Garding inequality says that they converge in all Sobolev
$L^{2}_{k}$-norms (with respect to the $g_{0}$-metric on $X$), and
hence uniformly on compact subsets $K\subset X_{+}\cup (0,\eps ]\times
Y\subset X$, in particular. Now, there is a diffeomorphism
$\phi_{i}:X(R_{i})\rarr X$, which takes the piece $X_{+}\cup (0,R_{i}]\times
Y$ to the piece $X_{+}\cup (0,\eps ]\times Y$, taking the metric
$g_{R_{i}}:=g_{i}$ on $X(R_{i})$ to the metric $g_{i}$ on $X$. Note
that $\phi_{i}$ are identity on $X_{+}$, and just a scaling of
$\frac{\eps}{R_{i}}$ of the $t$ variable on $(0, R_{i}]\times Y$ and 
identity on the $Y$ variable. So, on the piece $(0,\eps]\times Y$, the
1-form $ds$ has length $\frac{\eps}{R_{i}}$ with respect to $g_{i}$,
but length 1 with respect to $g_{0}$. The 1-forms on $Y$ have the same
length with respect to $g_{i}$ and $g_{0}$. It is now more
convenient to change the variable $t$ to $R-t$, which replaces 
 $X_{+}\cup (0, R]\times Y$ with the isometric manifold $X_{+}\cup [0,
R)\times Y$, where 
$\partial X_{+}$ is glued to $\{0\}\times Y$, and similarly $s$ to
$\eps - s$ taking $(0,\eps]$ to
$[0,\eps)$. Let us define $X_{i+}=X_{+}\cup [0, R_{i})\times Y$, and $\phi_{i}$
is the composite diffeomorphism $\phi_{i}:X_{i+}\rarr X_{+}\cup
[0,\eps)\times Y$. Also, since $\phi_{i}$ and $\phi_{j}$ are both 
identity on $Y$, and 
$\phi_{i}^{*}ds = \frac{\eps}{R_{i}}dt,\;\phi_{j}^{*}ds =
\frac{\eps}{R_{j }}dt$ for $j\geq i$, where $s\in [0,\eps)$, it follows that 
for $j\geq i$ we have, for any i-form $\omega$ on
$[0,\eps)\times Y$,  the inequalities:
\begin{eqnarray}
\norml\phi_{i}^{*}\omega\normr_{g_{i},X_{i+}}&\leq&
\norml\omega\normr_{g_{0},X_{+}\cup [0,\eps)\times Y}\nonumber\\
\norml\phi^{*}_{i}\omega - \phi_{j}^{*}\omega\normr_{g_{j},X_{j+}}&\leq & 
\eps\left(\frac{1}{R_{i}}-\frac{1}{R_{j}}\right)\norml
\omega\normr_{g_{0}, X_{+}\cup [0,\eps)\times Y}  
\label{normineq}
\end{eqnarray}

Note that all the 
$X_{i+}$ are isometrically embedded in the non-compact manifold
with infinite end $ X_{\infty}:=X_{+}\cup [0,\infty)\times Y$, where
$\{0\}\times Y$ is glued to $\partial X_{+}$,
with the 
 complete metric 
defined by $g_{+}$ on $X_{+}$ and $dt^{2}\times g_{Y}$ on the infinite tube 
$[0,\infty)\times Y$ (call this metric $g_{\infty}$), so that
$g_{\infty |X_{i}}= g_{i}$. 
Extending $\phi_{i}^{*}\omega_{i}$ on $X_{i+}$ by $0$ to all of
$X_{\infty}$ defines an $L^{2}(g_{\infty})$ form on $X_{\infty}$,
which we continue to denote by the same symbol. Now let 
$\til{\omega}_{i}$ be the $\Delta_{g_{\infty}}$-harmonic part of
$\phi_{i}^{*}\omega_{i}$. This is possible in view of the Kodaira
decomposition:
$$L^{2}={\cal H}^{2}_{\infty}\oplus \overline{d\Lam_{c}}\oplus 
\overline{\delta\Lam_{c}}$$
which is always true for a {\em complete} Riemannian metric (see [Ko]). (For
such a complete metric $g$, $\kernel{\Delta_{g}}$ is the same as the 
space of closed and co-closed forms.) Now
let $K$ be a compact subset of
$X_{i+}$, and hence $X_{\infty}$. Since for $j\geq i$, we have $X_{i+}\subset
X_{j+}\subset X_{\infty} $ and
$g_{j}=g_{i}$ on $X_{i+}$, we have the inequalities of sup
($C^{0}$) norms :
\begin{eqnarray*}
\norml\til{\omega}_{i}-\til{\omega}_{j}\normr_{\infty,K,g_{\infty}}
&\leq&\norml\phi_{i}^{*}\omega_{i}-\phi_{j}^{*}\omega_{j}
\normr_{\infty,K,g_{\infty}}\\
&=& \norml\phi_{i}^{*}\omega_{i}-\phi_{j}^{*}\omega_{j}
\normr_{\infty,K,g_{j}}\\
&\leq&\norml
\phi_{i}^{*}\omega_{i}-\phi_{j}^{*}\omega_{i}\normr_{\infty,K ,g_{j}}
+ \norml \phi_{j}^{*}\omega_{i} -
\phi_{j}^{*}\omega_{j}\normr_{\infty,K,g_{j}}\\ 
&\leq&
\eps\left(\frac{1}{R_{i}}-\frac{1}{R_{j}}\right)\norml\omega_{i}
\normr_{\infty,\phi_{j}(K),g_{0}} +\norml\omega_{i}-\omega_{j}
\normr_{\infty,\phi_{j}(K),g_{j}}\\
&\leq& \eps\left(\frac{1}{R_{i}}-\frac{1}{R_{j}}\right)\norml\omega_{i}
\normr_{\infty,\phi_{j}(K),g_{0}} + \norml\omega_{i}-\omega_{j}
\normr_{\infty,\phi_{j}(K),g_{0}}
\end{eqnarray*}
by using the two inequalities ~(\ref{normineq}) above.
This shows that
the $\Delta_{g_{\infty}}$-harmonic forms $\til{\omega}_{i}$ are uniformly
Cauchy on compact subsets of $X_{\infty}$, and hence converge
uniformly on compact sets to some $\Delta_{g_{\infty}}$-harmonic form
$\til{\omega}$. Also, 
\begin{eqnarray*}
\norml\til{\omega}\normr_{g_{\infty},X_{\infty}}&=&
\lim_{i\rarr\infty}\norml\til{\omega}_{i}\normr_{g_{i},X_{i}}\\&\leq& 
\lim_{i\rarr\infty}\norml\phi_{i}^{*}\omega_{i}\normr_{g_{i},X_{i}}\\
&\leq&\lim_{i\rarr\infty}\norml\omega_{i}\normr_{g_{0},X_{+}\cup
[0,\eps)\times Y} < \infty 
\end{eqnarray*}
shows that $\til{\omega}$ is an $L^{2}$ 2-form with respect to
$g_{\infty}$ on $X_{\infty}$, so it is in ${\cal H}^{2}_{\infty}$. 
However, the Kodaira decomposition above shows that this space is
precisely the image of $H^{2}_{c}(X_{\infty})$ in $H^{2}(X_{\infty})$.
However, from the definition of $X_{+}$, 
this  is precisely the image of $H_{2}(X_{+},\partial X_{+})$
in $H^{2}(X_{+})$, which is the kernel of the restriction map  
$H^{2}(\til{\Sigma}\times D^{2})~\rarr~H^{2}(\til{\Sigma}\times
S^{1})$, which is zero (e.g. by Kunneth formula). Thus $\til{\omega}=0$.

Now, since $\nu_{\eps}(\til{\Sigma})$ is a trivial bundle, we may find
a copy of $\til{\Sigma}$, call it $\til{\Sigma}_{1}$ which lies in 
$\partial X_{+}=\partial  \nu_{\eps}(\til{\Sigma})$, and is homologous
to $\til{\Sigma}$ in $X$. Note that
$\til{\Sigma}_{1}$ therefore lies in  $X_{i+}$ for all $i$, and 
in $X_{\infty}$,
satisfying $\phi_{i}^{*}(\til{\Sigma}_{1})=\til{\Sigma}_{1}$ for all
$i$. Now,
\begin{eqnarray*}
\lim_{i\rarr\infty}\left(\omega_{i}\cup (dH-E)\right)&=& \lim_{i\rarr\infty}\int_{\;\til{\Sigma}}\omega_{i} = 
\lim_{i\rarr\infty}\int_{\;\til{\Sigma_{1}}}\omega_{i}\\
&=&\lim_{i\rarr\infty}\int_{\;\til{\Sigma}_{1}}\phi_{i}^{*}\omega_{i}=
\lim_{i\rarr\infty}\int_{\;\til{\Sigma}_{1}}\til{\omega}_{i}\\
&=&\int_{\;\til{\Sigma}_{1}}\til{\omega}=0
\end{eqnarray*}

Consequently, since $c_{1}(L)=3H-E$, we have
\begin{eqnarray*}
\lim_{i\rarr\infty}\left(\omega_{i}\cup c_{1}(L)\right)&=&
\lim_{i\rarr\infty}\omega_{i}\cup (dH - E) - (d-3)(\omega_{i}\cup H)\\
&=& 0 - (d-3)  < 0
\end{eqnarray*}
since $d > 3$ by assumption. This proves the proposition.\hfill $\Box$

\begin{corollary}{\rm There exists a metric $g$ on $X$, which is a
product in a tubular neighbourhood of $Y$, such that $\#
M_{c,\delta}(g_{R})\neq 0$ for $R$ large enough. }
\label{nonempty2}
\end{corollary}

\vspace{.25in}
\noin
{\bf Proof:} Since $Y$ has a product neighbourhood, one can put a
product metric on an $\eps$-neighbourhood of $Y$, and extend it to all
of $X$ by using a partition of unity. The rest follows from the Corollary
~\ref{nonempty} and the proposition above.\hfill $\Box$

\subsection{Temporal Gauge Solutions}
We go back to the setting of the tubing construction of
Def.~\ref{tubing}, and do some analysis on the tube portion $[-R,R] \times
Y$. Look at the restrictions of the $U(2)$ bundles $W_{\pm}$ 
coming from the chosen $\mbox{Spin}_{c}$ structure on $X$  to
$[-R,R]\times Y$, viz. $W_{\pm|[-R,R]\times Y}$. These are both isomorphic via
$c(dt)$, 
Clifford multiplication by the (unit length) 1-form $dt$, and hence
may both be regarded as pullbacks via the projection map $[-R,R]\times
Y~\rarr~Y$ of a $U(2)$-bundle $W_{3}$ on $Y$.  If $A$ is a connection
on $[-R,R] \times Y$, then $A_{|\{t\}\times Y}:=A(t)$ is a connection for
$L_{|Y}$ for each $t$, and maybe regarded as a path in ${\cal A}_{Y}$,
the affine space of $U(1)$ connections on $L_{|Y}$. Similarly, if
$\Phi\in \Gamma (W_{+})$, one may regard the restriction to the
$t$-slice $\Phi(t)$ as a path of sections of $\Gamma(W_{3})$. Note
that the metric on $[-R,R]\times Y$ is a product, and hence induces
the same metric on each slice $\{t\}\times Y$. 
\begin{definition}{\rm We say a connection $A$ on $\R\times Y$ is in
{\em temporal gauge} if it has no $dt$ component. (We are of course fixing a
reference connection $A_{0}$, as always). Say it is {\em translation
invariant in a temporal gauge} if it is the pullback of a connection
on $L_{|Y}$. Similar definitions make sense for $[0,R)\times Y$, and 
$[-R,R]\times Y$. Finally, a solution $(A,\Phi)$ to the Seiberg-Witten
equations on $(X_{R}=X,\; g_{R})$ will be called a {\em temporal gauge
solution} if $A$ is in 
temporal gauge on the tube $T=[-R,R]\times Y$, and similarly for
translation invariant temporal gauge solutions}
\label{temporaldef}
\end{definition}

Clearly, if a connection $A$ on $L_{|\R\times Y}$ is in temporal
gauge, it can be recovered from the path $A(t)$, since it has no $dt$
component. 
\begin{remark}{\rm If $A=A_{R}$ is any connection on $X_{R}$ (see
~\ref{tubing} ), there exists a gauge transformation $g(R)$ in the connected
component of $id_{X}$ in the gauge group ${\cal G}$ such that $g(R)A$ is in
temporal gauge on $[-R, R]\times Y$. Further, if $A_{R|X_{+}}$  is
a fixed connection $A_{+}$ on $L_{|X_{+}}$ independent of $R$, 
we may choose $g(R)$ to be the identity map $id_{X_{+}}$ on the piece
$X_{+}$ for all $R$. Similarly for $X_{-}$. }
\label{temporalremark}
\end{remark}

\vspace{.25in}
\noin
{\bf Proof of remark:} 

Let $A_{0}(y,t)dt$ be the $dt$ component of $A$ on 
$[-R,R] \times Y$. There is a $C^{\infty}$ function $h_{R}(y,t)$ on 
$[-R,R] \times Y$ which satisfies : 
$A_{0}(t,y)=\frac{\partial h_{R}(t,y)}{\partial t}$.
Choose a $C^{\infty}$-function $f$ extending $h_{R}$ to $X=X_{R}$. If
$A_{R} (R, y)= A_{+}(R ,y)$ is independent of $R$, we may choose
$h_{R}(t,y)$ to be such that $h(R,y)=0$ for all $R$, and choose 
$f$ to be identically 0 on $X_{+}$ for all $R$. Now take
the gauge transformation $g(R)=\exp (2\pi i)f(t,y)$. Now, since
$$g(R)^{-1}dg(R)~=~2\pi i \left(
\frac{\partial h_{R}(t,y)}{\partial t}\right)dt=2\pi iA_{0}(y,t)dt$$ 
on $[-R,R]\times Y$, it follows that $gA=A-\frac{1}{2\pi i} g^{-1}dg$
has no $dt$ component on $[-R,R] \times Y$. \hfill $\Box$

Let $\partial_{A(t)}:\Gamma(W_{3})\rarr \Gamma(W_{3})$ be the induced
Dirac operator on $W_{3}$, with respect to the connection $A(t)$. In
view of Remark~\ref{temporalremark} above, any solution $(A,\Phi)$ to
the Seiberg-Witten equations on $(X_{R}, g_{R})$ can be assumed to
be a temporal gauge representative in its gauge equivalence class.  
Our first goal is to consider the
restriction of such a temporal gauge solution to the tube  
$T=[-R, R]\times Y$, and view it as a time-dependent solution
$(A(t),\Phi(t))$ of some equations on $Y$ involving $\partial_{A(t)}$ etc.

\begin{lemma}{\rm Let us denote by $*_{Y}$ the star operator on $Y$ defined by
the metric $g_{Y}$ (induced from $g_{R}$ on $X$ as in ~\ref{tubing}),
and let a 2-form $\omega$ on $X$ be expressed as $\omega = dt\wedge
\phi + \psi$ on the tube $T=[-R,R]\times Y$ where $\phi,\,\psi$ are
devoid of $dt$. Then $$*\omega = *_{Y}\phi + dt\wedge *_{Y}\psi$$
where $*$ is the star operator of $g_{R}$ on $X$. }
\label{starY}
\end{lemma}

\vspace{.25in}
\noin
{\bf Proof:} Is a straightforward exercise, since $g_{R}=dt^{2}\times
g_{Y}$ and so $dV_{g_{R}}=dt dV_{Y}$, and the star operator is
characterised by the diagram in the subsection  1.3 (Self-Duality)
of \S 1. \hfill $\Box$

\begin{corollary}{\rm Let $A$ be a connection on $X=X_{R}$, and let
$F_{A}$ be its curvature. Let $T$ denote the tube $[-R, R]\times Y$, with the
product metric $g_{R}=dt^{2}\times g_{Y}$. Then, 
if $A$ is translation invariant in a temporal gauge,  we have the equality of 
pointwise norms :
$$(F_{A}^{+}(y), F_{A}^{+}(y))_{y,g_{R}}=
(F_{A}^{-}(y), F_{A}^{-}(y))_{y,g_{R}}$$.
}
\label{plusminus}
\end{corollary}

\vspace{.25in}
\noin
{\bf Proof:} Write $F_{A} = dt\wedge \frac{dA}{dt} + F_{A,Y}$. 
By the Lemma above, 
\begin{eqnarray*}
F_{A}^{+}=\frac{1}{2}\left(F_{A}+*F_{A}\right)
&=&\frac{1}{2}\left[dt\wedge\left(
\frac{dA}{dt}+ *_{Y}F_{A,Y}\right) + \left(F_{A,Y} + *_{Y}\frac{dA}{dt}\right)
\right]\\
F_{A}^{-}=\frac{1}{2}\left(F_{A}+*F_{A}\right)&=&
\frac{1}{2}\left[dt\wedge\left(
\frac{dA}{dt}-*_{Y}F_{A,Y}\right) + \left(F_{A,Y} -*_{Y}\frac{dA}{dt}\right)
\right]
\end{eqnarray*}
Thus, since $\frac{dA}{dt}=0$ for $A$ translation-invariant and in
temporal gauge, we have the result.
\hfill $\Box$

\vspace{.25in}
The proof above also shows that the spaces of $*$-self dual and
antiself dual 2-forms on $T$ which are $t$-translation invariant are
both isomorphic to $\Omega^{1}(Y)$. The self-dual one is given as
$\frac{1}{2}(dt\wedge *_{Y}\omega (y) + \omega (y))$ and the
anti-selfdual one as  $\frac{1}{2}(-dt\wedge *_{Y}\omega (y) + \omega (y))$.

For notational convenience, we denote the isomorphism $\omega\mapsto
\frac{1}{2}(-dt\wedge *_{Y}\omega + \omega )$ by 
$\theta~:~\Omega^{1}(Y)\rarr~\Omega^{2+}_{\mbox{inv}}(T)$, where the right 
hand side is the space of translation invariant self-dual forms on $T$. 

The Clifford structure map $\gamma:\Lambda^{1}(T)\otimes \C
\rarr \mbox{Hom}(W_{+},W_{-})$ restricts to the Clifford isometry 
$\til{\gamma}:\Lambda^{1}(Y)\otimes\C\rarr\mbox{End}^{0}(W_{3})$. It is 
easily checked that $\rho\circ \theta = \til{\gamma}$. 
Finally, we recall the pairing $\sigma$ defined by
$\sigma(\Phi,\Psi)=i\left(h_{+}(-,\Psi)\Phi-\frac{1}{2}
h_{+}(\Phi,\Psi)\mbox{Id}\right)$
over the tube $T$:
$$\sigma :W_{+}\otimes\overline{W}_{+}\rarr
\mbox{End}^{0}(W_{+})\simeq \Lam^{2+}\otimes \C$$
 Since
$W_{+}=\pi^{*}(W_{3})$, we have the pairing:
$$\tau : W_{3}\times \overline{W}_{3}\rarr \Lambda^{1}(Y)\otimes \C$$  
where $\til{\gamma}\circ\tau = \sigma$. 
So, if we regard a section $\Phi\in \Gamma(W_{+|T})$ as a section
in $\Gamma(\pi^{*}W_{3})$, which is the same as a path of
sections $\Phi(t)$ in $\Gamma(W_{3})$, then by definition, 
$$\sigma(\Phi(t,y),\Phi(t,y))=\til{\gamma}(\tau(\Phi(t)(y),\Phi(t)(y)))$$
 
Now we look at the Dirac operator on $T$. The covariant derivative
with respect to $A$ and compatible with the Levi-Civita
connection of $g_{R}$ is given by :
$$\nabla_{A,dt}= i(A\vdash dt)\otimes(-)\;+\; \frac{\partial}{\partial t}$$
where $\vdash$ denotes contraction. Since $A$ is in temporal
gauge, $A\vdash dt =0 $, and so
$\nabla_{A,dt}=\frac{\partial}{\partial t}$. Of course, viewing
$\Phi(t,-)$ as a path of sections $\Phi(t)$ of $W_{3}$, we have 
$$\frac{\partial\Phi(t,-)}{\partial t}=\frac{d\Phi(t)}{dt}$$
Also, the Clifford
multiplication $c(dt)$ is what makes $W_{+}$ isomorphic to
$W_{-}\simeq \pi^{*}(W_{3})$. So the Dirac
operator for the bundle $W$ on the tube $T$ with induced
$\mbox{Spin}_{c}$ structure on $T$ reads as :
$$D_{A} =\frac{d}{dt}+\partial_{A(t),Y}$$
where $\partial_{A(t),Y}$ is the Dirac operator for the induced 
$\mbox{Spin}_{c}$ structure on $Y$. 
It follows that time-dependent Seiberg-Witten equations read on $Y$ as
(assuming,
as usual that $A$ is in temporal gauge):
\begin{eqnarray}
\left(\frac{dA}{dt} + *_{Y}F_{A,Y}\right)-\tau(\Phi,\Phi)-
*_{Y}\frac{\delta}{2} =0 \nonumber\\
\frac{d\Phi(t)}{dt}= -\partial_{A(t),Y}\Phi(t)
\label{tdependent}
\end{eqnarray}
where $\delta\in \Omega^{2}(Y)$.
Hence we have the :
\begin{proposition}{\rm If $(A,\Phi)$ is a temporal gauge 
solution to the Seiberg-Witten 
equations on the tube $T=\R\times Y$ or $[-R, R]\times Y$, then 
$(A(t),\Phi(t))$ is 
a path in ${\cal A}_{Y}\times \Gamma(W_{3})$ which is a trajectory of
the equations ~(\ref{tdependent}) 
\label{trajectories}
}
\end{proposition}
We now note that the equations ~\ref{trajectories} above are the
gradient-flow equations of a functional defined on 
${\cal A}_{Y}\times\Gamma(W_{3})$. 

\begin{proposition}{\rm
The equations ~(\ref{tdependent}) are the gradient
flow equations for the functional
defined on ${\cal A}_{Y}\times \Gamma(W_{3})$ by 
$$C_{\delta}(A,\Phi)=-\frac{1}{2}\left(\int_{\,Y}(A-B)\wedge (F_{A} -\delta) +
\int_{\,Y}\left(\Phi,\partial_{A}\Phi\right)_{W_{3},y} dV_{Y}\right)$$
where $B$ is a reference connection on $L_{|Y}$ and the integrand on
the extreme right is the  inner product on the fibre $W_{3,y}$, and
$\delta$ is a 2-form on $Y$. Hence, $(A(t),\Phi(t))$ satisfying
~(\ref{tdependent}) is the gradient flow for this functional, and
$C_{\delta}$ is  monotonically increasing along this trajectory. 
}
\label{functional}
\end{proposition}

\vspace{.25in}
\noin
{\bf Proof:}  For simplicity, denote $A-B$ as $A$, where $B$ is the reference
connection.
Let $\{e_{i}\}_{i=1}^{2}$  be a $(-,-)_{W_{3}}$ unitary frame for $W_{3}$. We 
recall that $\sigma(\Phi,\Psi)=i\left(h_{+}(-,\Psi)\Phi-\frac{1}{2}
h_{+}(\Phi,\Psi)\mbox{Id}\right)$.  
The foregoing definitions lead to the following 
identity (since $\til{\gamma}$ maps into traceless endos):
\begin{eqnarray*}
<\tau(\Phi,\Phi),\omega>&=&<\sigma(\Phi,\Phi),\til{\gamma}(\omega)>=
\frac{1}{2}\mbox{Tr}((-,\Phi)\Phi\circ\til{\gamma}(\omega)^{t})\\
&=&\frac{i}{2}\sum_{j=1}^{2}((\til{\gamma}(\omega)^{t}e_{j},\Phi)\Phi,e_{j})=
\frac{i}{2}(\Phi,\til{\gamma}(\omega)\Phi)
\end{eqnarray*}
where the round brackets denote  the hermitian inner product $(-,-)_{W_{3},y}$
on the fibre $W_{3,y}$ and the angular brackets denote the Riemannian 
inner product on $Y$. In view of the calculation above, and that 
$\partial_{A} = \sum_{j=1}^{3}\til{\gamma}(\omega_{j})\nabla_{A,\omega_{j}}$
the inner product $(\Phi,\;\partial_{A}\Phi)_{W_{3}}$ 
 satisfies:
\begin{eqnarray*}
(\Phi,\;\partial_{A}\Phi)_{W_{3}}-((\Phi,\;\partial_{B}\Phi)_{W_{3}}&=&
-i\sum_{j=i}^{3}(\Phi, <A,\omega_{j}>\til{\gamma}(\omega_{j})\Phi)_{W_{3}}\\
&=& -i\sum_{j=1}^{3}<A, \omega_{j}>(\Phi,\til{\gamma}(\omega_{j}))_{W_{3}}
=-2<\tau(\Phi,\Phi),A>
\end{eqnarray*}
where $\{\omega_{j}\}_{j=1}^{3}$ is a local orthonormal frame for $\Lambda^{1}(Y)$.

 Thus the integrand of
$C_{\delta}$ becomes :
$$Q(A,\Phi)=-\left[\frac{1}{2}<A, *_{Y}(F_{A}-\delta)> + 
\frac{1}{2}(\Phi, \partial_{B}\Phi)_{W_{3},y}-<\tau(\Phi ,\Phi), A>\right]$$
Thus,
\begin{eqnarray*}
\frac{\partial Q}{\partial A}&=&-\left(<-,*_{Y}(F_{A}-\frac{1}{2}\delta)> -
 <-, \tau(\Phi,\Phi)>\right)\\
\frac{\partial Q}{\partial\Phi} &=&
-\frac{1}{2}\left((-,\partial_{A}\Phi)_{W_{3},y}+ 
(\Phi, \partial_{A}-)_{W_{3},y}\right)\\
&=&  -(-,\partial_{A}\Phi)_{W_{3},y}
\end{eqnarray*}
by using Stokes formula, $F_{A}=dA$ and the self-adjointness of $\partial_{A}$ on $Y$.
Thus the gradient flow equations for the functional, viz.,
\begin{eqnarray*}
\frac{dA}{dt} =  \frac{\partial Q}{\partial A} \\
\frac{d\Phi}{dt}= \frac{\partial Q}{\partial\Phi}
\end{eqnarray*}
lead, respectively, to the required equations~(\ref{tdependent}), and the proposition is
proved.\hfill $\Box$. 

We next investigate what happens to $C_{\delta}$ under the Gauge group
action. We shall drop the subscript $\delta$ from $C_{\delta}$ for
notational convenience.

\begin{proposition}{\rm Under a gauge transformation $g\in
\mbox{Map}(Y, S^{1})$, we have the 
transformation formula:
$$C(gA,g\Phi)=
C(A,\Phi)+2\pi^{2}[g]\cup\left(c_{1}(L)-\frac{1}{2\pi}[\delta_{\cal H}]\right)$$
where $[g]$ denotes the cohomology class of $g$ in $H^{1}(Y,\R)$.
Since $[g]$ is an integral cohomology class, this shows 
that on a gauge orbit, $C$ is well defined in $\R / 2\pi^{2}\Z \cup 
\left(c_{1}(L)-\frac{1}{2\pi}[\delta_{\cal H}]\right)$. 
}
\label{cgaugelaw}
\end{proposition}

\vspace{.25in}
\noin
{\bf Proof:}  We apply the formulas $gA=
A-g^{-1}dg,\;\partial_{gA}(g\Phi)=g\partial_{A}\Phi$, and $F_{gA}=F_{A}$
to compute:
\begin{eqnarray*}
C(gA, g\Phi) &=& C(A,\Phi) +\frac{1}{2}\int_{\,Y} g^{-1}dg\wedge (F_{A}
-\delta)\\
&=& C(A,\Phi)  +\pi[g]\cup 2\pi(c_{1}(L)-\frac{1}{2\pi}[\delta_{\cal H}])\\
&=& C(A,\Phi) +2\pi^{2}[g]\cup (c_{1}(L)-\frac{1}{2\pi}[\delta_{\cal H}])
\end{eqnarray*}
proving the proposition. \hfill $\Box$

Now we are ready for the main proposition of this section.

\begin{proposition}{\rm In the notation of Definition~\ref{tubing},
assume that $M_{c,\delta}(g_{R})$ is non-empty for $R$ large enough.
Then there exists a solution $(A,\Phi)$ to the Seiberg-Witten
equations on $\R\times Y$ which is translation invariant in a temporal
gauge. }
\label{translationinvariant}
\end{proposition}

\vspace{.25in}
\noin
{\bf Proof:} Let $(A_{R},\Phi_{R})$ be a solution to the
Seiberg-Witten equations on $(X_{R}, g_{R})$ for sufficiently large
$R$. By gauge transforming if necessary (see Remark~\ref{temporalremark}), let us assume that all the $A_{R}$ are in temporal gauge.
Take  reference connections $B_{R}$ on the bundle $L\rarr X_{R}$ such
that the restrictions $B_{R|X_{\pm}}:=B_{\pm}$ are fixed, independent
of $R$. 

Let us denote the change in the functional $C$ along the tube 
$T=[-R, R]\times Y \subset X_{R}$ by :
$$l_{A,\Phi}(R) = C(A_{R}(R),\Phi_{R}(R))-C(A_{R}(-R),\Phi_{R}(-R))$$
By the Proposition~\ref{functional}, we
have that $C(A_{R}(t),\Phi_{R}(t))$ are monotonic functions of $t$. 
 Now, $(A_{R},\Phi_{R})$ will
restrict to solutions on $X^{\pm}$. By construction, the scalar
curvature of $g_{R}$ has the same
infimum on $X_{R}$ for all $R$. Hence, by \S 5.2 of [PP], there is a uniform
$C^{0}$-bound for all  $\Phi_{R}$ on $X_{R}$ 
independent of $R$. The (compactness)
argument of \S 5.2 of [PP]
shows that there exist gauge transformations $h_{R}^{\pm}$ of
$X_{\pm}$ such that :
$$h_{R}^{\pm}A_{R} - B_{R}^{\pm}=h_{R}^{\pm}A_{R} - B^{\pm}$$
are both bounded in Sobolev $L^{2}_{k}$-norm (for $k$ suitably large) 
uniformly for all $R$. By Sobolev's Lemma, this implies uniform $C^{0}$ bounds
on zeroth and first derivatives of $A_{R},\Phi_{R}$ for all $R$. 
Since $C(h_{R}^{\pm}A_{R}(\pm R),\Phi_{R}(\pm R))$ 
are the evaluations of $C$ (which only involves derivatives upto first
order) to the ends $\{\pm R\}\times
Y$ of the tube $[-R, R]\times Y$, we have uniform bounds for  
$C(h_{R}^{\pm}A_{R}(\pm R),\Phi_{R}(\pm R))$
independent of $R$. Now let $\gamma$ be a 1-cycle Poincare-dual to 
$c_{1}(L)-\frac{1}{2\pi}[\delta_{\cal H}]$ in $Y$. It is the
intersection with $Y$ of the 2-cycle $\Gamma$ which is Poincare-dual
to  $c_{1}(L)-\frac{1}{2\pi}[\delta_{{\cal H}}]$ in $X$. 

Let $i_{\pm}:Y\hookrightarrow X_{\pm}$ denote the inclusions (as ends),
and $[h_{R}^{\pm}]$ denote the cohomology classes of $h_{R}^{\pm}$ in
$H^{1}(Y)$, or $H^{1}(X_{\pm})$. Then,
\begin{eqnarray*}
<\left(c_{1}(L)-\frac{1}{2\pi}[\delta_{\cal H}]\right)\cup [h_{R}^{\pm}],[Y] >
&=&<i_{\pm}^{*}\left(c_{1}(L)-\frac{1}{2\pi}[\delta_{\cal H}]
\right)\cup [h_{R}^{\pm}],[Y]> \\
=<i_{\pm}^{*}\left(c_{1}(L)-\frac{1}{2\pi}
[\delta_{\cal H}]\right)\cup [h_{R}^{\pm}], \partial_{\pm}[X_{\pm}]>
&=& <\delta_{\pm}i_{\pm}^{*}\left(c_{1}(L)-\frac{1}{2\pi}[\delta_{\cal
H}]\right)\cup [h_{R}^{\pm}],[X_{\pm}]>=0
\end{eqnarray*}
where $[Y],\,[X_{\pm}]$ denote orientation classes, and
$\partial_{\pm}:H_{4}(X_{\pm},Y)\rarr H_{3}(Y)$ and
$\delta_{\pm}:H^{3}(Y)\rarr H^{4}(X_{\pm},Y)$ denote respectively the
connecting homomorphisms in the long exact homology and cohomology
sequences of the pair $(X_{\pm},Y)$. Thus, by
Proposition~\ref{cgaugelaw} 
we have $$C(h_{R}^{\pm}A_{R}(\pm R), h_{R}^{\pm}\Phi_{R}(\pm R))= 
C(A_{R}(\pm R), \Phi_{R}(\pm R))$$
Thus we have a uniform bound $M$ on $l_{A_{R},\Phi_{R}}$
independent of $R$. 

 Now let $R$ be a positive integer,
say $R=N$, and denote by $\Delta_{i}$ the change in $C$ across 
$[i-1,i]\times Y$, viz. 
$$\Delta_{i}= C(A_{N}(i),\Phi_{N}(i)) - C(A_{N}(i-1),\Phi_{N}(i-1))$$
We saw in Proposition~\ref{trajectories} that $(A_{N},\Phi_{N})$, being
 solutions to Seiberg-Witten equations and $A_{N}$  being in temporal
gauge implied that they were (time-dependent) solutions to the
equations ~(\ref{tdependent}), and  the Proposition~\ref{functional}
then implied that $C$ was  monotonic increasing in time for these solutions. 
 Thus all the $\Delta_{i}$ are non-negative. Let
$\Delta_{min,N}=\min_{i}\Delta_{i}$. Hence we have :
$$2N\Delta_{\min,N}\leq \sum_{i=-N}^{+N}\Delta_{i}
=C(A_{N}(N),\Phi_{N}(N)) -
C(A_{N}(-N),\Phi_{N}(-N))=l_{A_{N},\Phi_{N}}\leq M$$
for all $N$. 

Hence $\lim_{N\rarr\infty}\Delta_{\min ,N}= 0$. Denote by 
$(A_{(N)},\Phi_{(N)})$ the restriction of $(A_{N},\Phi_{N})$ to the interval 
$[i-1, i]\times Y$ on which 
$\Delta_{\min ,N}=\left|\Delta_{i}\right|$. 
This may be viewed as a solution on $[0,1]\times Y$, denoted by
the same symbol  $(A_{(N)},\Phi_{(N)})$. As we saw above, we have

$$C(A_{(N)}(1),\Phi_{(N)}(1))-C(A_{(N)}(0),\Phi_{(N)}(0))
\leq \frac{M}{N}$$
which goes to $0$ as $N\rarr\infty$. The uniform $C^{0}$ bound
on $(A_{N},\Phi_{N})$, and hence $(A_{(N)},\Phi_{(N)})$ gives a solution 
(on passing to a
subsequence) $(A,\Phi)$ on $[0,1]\times Y$ for which 
$$\left(C(A(1),\Phi(1))-C(A(0),\Phi(0))\right)
=\lim \left(C(A_{(N)}(1),\Phi_{(N)}(1))-C(A_{(N)}(0),\Phi_{(N)}(0))
\right)=0$$
The monotonicity of $C$ across $[0,1]\times Y$ implies that
$(A,\Phi)$ is constant along $[0,1]\times Y$. This solution is
clearly therefore a translation invariant solution on
$[0,1]\times Y$, which extends to all of $\R\times Y$ by
time-translating for all times. \hfill $\Box$

\section{Proof of Thom's Conjecture}
We first need a lemma :
\begin{lemma}{\rm Let $Y=S^{1}\times\til{\Sigma}$, where
$\til{\Sigma}$ is a Riemannian 2-manifold of constant scalar
curvature $s$ and
genus $g \geq 1$. Assume the metric on
$\til{\Sigma}$ is normalised so that its volume is 1 (and thus 
$s=2\pi\chi(\til{\Sigma})=2\pi(2-2g)$). Let $Y$ have a metric $g_{Y}$
extending this metric on $\til{\Sigma}$, and let the infinite
tube $T=\R\times Y$ have a product metric $dt^{2}\times g_{Y}$,
and $L$ be the line bundle associated with a compatible
$\mbox{Spin}_{c}$ structure. 
Suppose there is a solution to the
Seiberg-Witten equations on $(T, g_{Y})$ which is translation
invariant in a temporal gauge. Then:
$$\left|\frac{1}{2\pi}\int_{\til{\Sigma}} F_{A}\right|
 \leq 2g -2 $$
}
\label{sigmatilde}
\end{lemma}

\vspace{.25in}
\noin
{\bf Proof:}  From the $C^{0}$ bound (see \S 5.2,
[PP]), the sup norm satisfies
$|\Phi|^{2}_{\infty}\leq 2\pi(2g-2) + \norml\delta\normr_{\infty}$.  Since
$\left|\sigma(\Phi,\Phi)\right|^{2}=\frac{1}{2}\left|\Phi\right|^{4}$,
we have 
$$|\sigma(\Phi,\Phi)|_{\infty}\leq
\frac{1}{\sqrt{2}}(\left(2\pi(2g -2)\right)+\norml\delta\normr_{\infty})$$

Thus, from the Seiberg-Witten equations:
$$\norml F_{A}^{+}\normr_{\infty ,Y}\leq
\frac{1}{\sqrt{2}}\left(2\pi(2g -2)\right) + 2\norml \delta\normr_{\infty} $$
However, by  Corollary~\ref{plusminus}, we have the pointwise
norm equality  $\norml F_{A}^{+}\normr =\norml F_{A}^{-}\normr$
because our solution is translation invariant in a temporal
gauge. Therefore,
$$\norml F_{A}\normr \leq \sqrt{2}\norml F_{A}^{+}\normr \leq 
2\pi(2g - 2) + O(\norml \delta\normr_{\infty}) $$
Thus 

$$\left|\frac{1}{2\pi}\int_{\,\til{\Sigma}}F_{A}
\right| \leq \frac{1}{2\pi}\left[ 2\pi(2g-2) +
O(\int_{\,\til{\Sigma}}\norml\delta\normr_{\infty})\right]$$

  This proves the lemma, since $\delta$ is arbitarily small.  
\hfill $\Box$

Now we can prove the main theorem.

\begin{theorem}[Kronheimer-Mrowka]{\rm If $\Sigma$ is an
oriented $2$-manifold smoothly embedded in $\C\Proj^{2}$ so as
to represent an algebraic curve of degree $d$, then the genus
$g(\Sigma)$ of $\Sigma$ satisfies:
$$g(\Sigma)\geq \frac{(d-1)(d-2)}{2}$$
}
\end{theorem}

\vspace{.25in}
\noin
{\bf Proof:} The cases of $d=1, 2$ are trivial, and $d=3$ is due
to Kervaire-Milnor (see reference [6] in [KM]). So we will
assume $d > 3$ in the sequel.   By Proposition ~\ref{negativecup}, there
exists a metric $g_{R}$ on
$X=\C\Proj^{2}\# d^{2}\overline{\C\Proj}^{2}$ such that $c_{1}(L)\cup
[\omega_{g_{R}}] < 0$. By Corollary ~\ref{nonempty2}, the moduli space 
$M_{c,\delta}(g_{R})\neq \phi $. (We just need to ensure that
the reference metric we started with on $X$ is a product metric
in a tubular neighbourhood of $Y=S^{1}\times\til{\Sigma}$.) By 
Proposition ~\ref{translationinvariant}, there is a solution on 
$\R\times Y$ which is translation invariant in a temporal gauge. By
the Lemma~\ref{sigmatilde} above,  
$\left| c_{1}(L).[\til{\Sigma}]\right|\leq 2g-2$, which implies
$c_{1}(L).\til{\Sigma}\geq 2-2g$. By the opening discussion of \S 
~\ref{peetwo} we have $c_{1}(L)=3H-E$, and by construction 
$[\til{\Sigma}]=dH-E$.
Thus $(3H-E).(dH-E)\geq 2-2g$. Applying
$H.E=0,\,H.H=1,\,E.E=-d^{2}$, we get $3d - d^{2}\geq 2- 2g$, i.e.
$g\geq \frac{(d-1)(d-2)}{2}$, proving the theorem.\hfill $\Box$.

\section{Appendix :Fredholm Theory}
\subsection{Preliminaries}
All Hilbert manifolds in the sequel are assumed to be second countable
and paracompact. We recall that a bounded operator 
$T:{\cal H}_{1}\rarr {\cal H}_{2}$ is said to be Fredholm if $\mbox{Ker}\,T$
and $\mbox{Coker}\,T$ are finite dimensional and the range
$\mbox{Im}\,T$ is closed. 
\begin{definition}{\rm We call a smooth map :
$$f:{\cal M}\rarr {\cal N}$$
a {\em Fredholm map} if the derivative $Df(x):T_{x}{\cal M}\rarr
T_{x}{\cal N}$ is a Fredholm operator for each $x\in {\cal M}$. }
\label{Fredmap}
\end{definition}
It is necessary to extend results like the implicit function
theorem in the finite dimensional case to the infinite
dimensional case, so as to make manifolds out of inverse images of
regular values etc. The key to doing it is the following proposition,
which enables one to construct a ``standard local model'' of a smooth
map whose derivative is given to be Fredholm at a point $p$, in
 a neighbourhood of $p$. Its main utility is to decompose a (non-linear)
smooth map with infinite-dimensional range into a {\em linear map} (with
infinite-dimensional range) and a {\em non-linear map} with {\em finite
dimensional} range, in a small neighbourhood of $p$.

\begin{proposition}{\rm Let ${\cal H}_{1}$ and ${\cal H}_{2}$ be
two Hilbert spaces, and let $f:{\cal H}_{1}\rarr{\cal H}_{2}$ be
a smooth map between them, such that $f(0)=0$. Assume that the derivative 
$T=Df(0)$ is a Fredholm operator from $T_{0}{\cal H}_{1}={\cal H}_{1}$
to $T_{0}{\cal H}_{2}={\cal H}_{2}$. {\em Note that we are only 
requiring the derivative to be Fredholm at a point, not that $f$
necessarily be a Fredholm map}. Then, with the orthogonal decomposition
${\cal H}_{1}=\mbox{Ker}\,T\oplus V_{1}$, ${\cal H}_{2}
=\mbox{Im}\,T\oplus V_{2}$, there exists a (non-linear) map 
$\phi~:~{\cal H}_{1}~\rarr~V_{2}=\coker{T}$ and a diffeomorphism $h:U\rarr h(U)$
for $U$ a neighbourhood of $0$ such that :
\begin{description}
\item[(i)] $f\circ h(n, v_{1}) = (Tv_{1}, \phi (n, v_{1}))$ for 
$n\in \mbox{Ker}\,T$, $v_{1}\in V_{1},\,(n,v_{1})\in U$. So $\phi$ is
smooth with finite dimensional range $V_{2}=\coker{T}$. In particular,
\item[(ii)] $f$ is a Fredholm map in the neighbourhood $U$ of $0$.
\item[(iii)] $\phi(0)=0$, $D\phi(0)=0$.
\item[(iii)] If $G$ is a group acting via an orthogonal linear action 
on ${\cal H}_{1}$ and ${\cal H}_{2}$, and $f$ is
$G$~-~equivariant, then $V_{1},\;V_{2}$ are $G$-invariant, and $h$
and $\phi$ are also $G$-equivariant.
\end{description}
}
\label{fredmodel}
\end{proposition}
\vspace{.25in}
\noindent
{\bf Proof:} 

Clearly if $f$ is $G$-equivariant, $T=Df(0)$ is
also $G$-equivariant, and the orthogonality of the action implies
that both $\mbox{Ker}\,T$ and $\im{T}$ being
$G$-invariant, their orthogonal complements $V_{1}$ and $V_{2}=\coker{T}$ are
$G$-invariant.

Note that $T_{|V_{1}}:V_{1}\rarr\mbox{Im}\,T$ is a bounded
bijective linear operator, so has a bounded inverse
$T^{-1}:\mbox{Im}\,T\rarr V_{1}$ by the open mapping theorem.
Let $\pi :{\cal H}_{2}\rarr\mbox{Im}\,T$ be the orthogonal
projection onto $\mbox{Im}\,T$. Let $\tilde{T}$ denote the composite
$T^{-1}\circ\pi:{\cal H}_{2}\rarr V_{1}$, and $\theta:{\cal H}_{1}\rarr
\kernel{T}$ be the orthogonal projection onto $\kernel{T}$. Consider
the map:$$
\begin{array}{clcl}
\chi :&{\cal H}_{1}&\rarr &{\cal H}_{1}\\
 &x&\mapsto &(\theta(x), \tilde{T}(f(x))
\end{array}
$$
Then, since $\theta(n)=n$ for $n\in \ker{T}$, we compute:
\begin{eqnarray*}
D\chi(0)(n,v)&=&(n, \tilde{T}\circ Df(0)(n,v)) =(n,\tilde{T}\circ T(n,v))\\
&=&(n, T^{-1}\circ\pi\circ T(n,v))= (n, T^{-1}\circ T(n,v))=(n,v)
\end{eqnarray*}
Note also that in the $G$-setting, $\chi$ is $G$-equivariant. 
Now $D\chi(0)=\mbox{id}$ implies by the (infinite-dimensional) inverse function
theorem that there exists a ball $V= B(0,\del)$ around the origin
in ${\cal H}_{1}$ on which $\chi$ is a diffeomorphism onto its
image. Let $U_{1}=B(0,\eps)\subset\chi(V)$, and for $y\in U_{1}$, define
$\phi_{1}:U_{1}\rarr V_{2}$ by $\phi_{1}=\pi_{V_{2}}\circ f\circ \chi^{-1}(y)$.
Note that again, in the $G$-setting, $\phi_{1}$ is
$G$-equivariant since it is a composite of $G$-equivariant maps.
Now for $y\in U_{1}=B(0,\eps)$, we have, using $T(\theta(u),v)=T(v)$
and the definitions above :
$$
\begin{array}{cl}
f\circ\chi^{-1}(y)&=(\pi\circ f\circ\chi^{-1}(y),\pi_{V_{2}}\circ f\circ
\chi^{-1}(y))=(T\circ T^{-1}\circ \pi\circ f\circ \chi^{-1}(y), \phi_{1}(y))\\
&= (T\circ \til{T}\circ f(\chi^{-1}(y),\phi_{1}(y))
=(T(\theta(\chi^{-1}(y)),\til{T} f(\chi^{-1}(y))), \phi_{1}(y))\\
&= (T\circ\chi\circ\chi^{-1}(y),\phi_{1}(y)) =(T(y),\phi_{1}(y))
\end{array}
$$
All we need to do now is extend $\phi_{1}:U_{1}\rarr V_{2}$ to 
$\phi :{\cal H}_{1}\rarr V_{2}$ in a $G$-equivariant manner, and
this is easily done by using the map :
\begin{eqnarray*}
\rho : &{\cal H}_{1} &\rarr {\cal H}_{1}\\
&x &\mapsto \frac{\eps x}{\eps + \psi(\norml x\normr )\norml x\normr}
\end{eqnarray*}
where $\psi:\R\rarr [0,1]$ is a $C^{\infty}$-map which is identically
zero for $|t| < \frac{\eps}{2}$ and identically = 1 for $|t| > \eps$.
Then $\rho$ is $G$-equivariant since $G$ preserves $\norml
\;\normr$; and maps ${\cal H}_{1}$ into $B(0,\eps)$ and is equal
to the identity map on $B(0,\frac{\eps}{2})$. Thus the map 
$\phi=\phi_{1}\circ \rho$ is a $G$-equivariant map agreeing with
$\phi_{1}$ on $B(0,\frac{\eps}{2})$. Now take 
$U=B(0,\frac{\eps}{2})$ and $h=\chi^{-1}$.
Clearly,  $\phi (0)= 0$ by construction, and since $Df(0)=T$, we have
$D\phi(0)=0$, and hence (i), (iii), and (iv) follow. To see
(ii) note that on $U$, we have $f$ equivalent to the map $(T, \phi)$
via the local diffeomorphism $h$ applied on the domain. Since $T$ is a 
linear isomorphism from $V_{1}$ to $\im{T}$, and $\phi$ is smooth with 
finite dimensional range, it is easy to check that $(T,\phi)$ is
Fredholm on  $U$. Thus $f$ is also Fredholm on $U$, proving (ii).\hfill $\Box$

\vspace{.25in}
\begin{corollary}{\rm Let $f,\,\phi, \,U,\, {\cal H}_{1},$ and
${\cal H}_{2}$ be as in the last proposition~\ref{fredmodel}.
Let $\til{U}:=h(U)$. Then, for $(a,\del )\in {\cal H}_{2}$, 
\begin{description}
\item[(i)] the germ of the inverse image of $(a,\delta)$, i.e. 
$f^{-1}(a,\del)\cap \til{U}$ is equivalent (i.e. via
an ambient diffeomorphism  $h:U\rarr\til{U}$) to 
$\phi_{b}^{-1}(\del)$, where 
\begin{eqnarray*}
\phi_{b} : &U_{b}\rarr &V_{2}\\
 &n\mapsto &\phi(n,b)
\end{eqnarray*}
$b:=T^{-1}(a)$ (uniquely defined as an element
of $V_{1}$), and $U_{b}=\theta(U \cap \pi_{V_{1}}^{-1}(b))$ is the
$b$-slice of $U$.
Thus, with the hypotheses on $f$ of the last proposition, a
local model for $f^{-1}(a,\del)$ near $0$ is given by the fibre 
of the {\em finite dimensional} smooth map $\phi_{b}$. In 
particular $f^{-1}(0)$ is locally homeomorphic to $\phi_{0}^{-1}(0)$.

\item[(ii)] For a fixed $a$ in $\im{T}\subset {\cal H}_{2},\;(a,\del)$
is a regular value for $f_{|\til{U}}$ if
and only if $\del$ is a regular value for the smooth map
$\phi_{b|U_{b}}$, {\em whose domain and range are  finite dimensional}. 
By the usual (finite-dimensional) Morse-Sard theorem applied to
$\phi_{b|U_{b}}$, it follows that for a fixed $a$, {\em there is a Baire
subset} of $\del\in V_{2}$ such that $(a,\delta)$ is a  regular
value for $f_{|\til{U}}$. In case $\del$ happens to be in
this Baire subset, the inverse image $f^{-1}(a,\del)\cap\til{U}$
is a smooth manifold of dimension = $\dim U_{b}\,-\dim\,V_{2}=
\dim\,(\kernel{T})-\dim\,(\coker{T})=\mbox{index}\,T$. ({\bf Note:} 
A {\em Baire} set is a countable intersection of open dense
sets.)

\item[(iii)] If $f:{\cal M}\rarr {\cal N}$ is a smooth Fredholm map between
Hilbert manifolds, the regular values of $f$ constitute a Baire
subset of ${\cal N}$. 
\end{description}
}
\label{transverse2}
\end{corollary}
\vspace{.25in}
\noin
{\bf Proof:} $f_{|\til{U}}$ has exactly the same properties as
$f\circ h_{|U}$ as far as regular values, inverse images of regular
values etc. are concerned. Thus ~\ref{fredmodel} implies (i) and
(ii). (iii) follows from second countability and the fact that 
countable intersections of Baire sets are Baire. \hfill $\Box$
\vspace{.25in}
\begin{corollary}{\rm Let everything be as in Proposition~\ref{fredmodel}.
If $f$ is a $G$-equivariant map, (with $G$ acting orthogonally
and linearly on both domain and range as before), then a local
(homeomorphic) model for $(\til{U}\cap f^{-1}(0))/G$ is 
$(U_{0}\cap\phi_{0}^{-1}(0))/G$ where $\phi_{0}=\phi(-, 0)$.}
\label{equivmodel1}
\end{corollary}
\vspace{.25in}
\noin
{\bf Proof:} By (i) of the corollary ~\ref{transverse2} above, 
$\til{U}\cap f^{-1}(0)$ is homeomorphic to 
$U_{0}\cap\phi_{0}^{-1}(0)$ via a $G$-equivariant local 
diffeomorphism $h:U\rarr\til{U}$ in the ambient space. Since
(iv) of ~\ref{fredmodel} implies that $\phi$ is $G$-equivariant,
so is $\phi_{0}=\phi(-,0)$, and the corollary follows. \hfill $\Box$

Thus we have proved the following :
\begin{proposition}(Local finite-dimensional model for a
neighbourhood of a singular point in the orbit space)
{\rm If $f:{\cal H}_{1}\rarr {\cal H}_{2}$ is 
a $G$-equivariant map between Hilbert spaces on which $G$ acts 
linearly and orthogonally, such that $T=Df(0)$ is a Fredholm
operator, then a local model for the germ of $f^{-1}(0)/G$ at
$0$ is the germ at $0$ of the finite dimensional object 
$\phi^{-1}_{0}(0)/G$, where $\phi_{0}:\kernel{T}\rarr\coker{T}$
is a smooth $G$-equivariant map between finite dimensional spaces
(with the restricted action of $G$). }
\label{equivmodel2}
\end{proposition}
\vspace{.25in}
Another lemma which will be useful in the sequel is the following:
\begin{lemma}{\rm 
Let $f:{\cal M}\rarr{\cal N}_{1}\times{\cal N}_{2}$ be a smooth map
between Hilbert manifolds, and let $f_{i}:=\pi_{i}\circ f$ where $\pi_{i}$ 
are the projection maps to ${\cal N}_{i}$ for $i=1,2$. Then 
\begin{description}
\item[(i)] $(a,b)\in {\cal N}_{1}\times {\cal N}_{2}$ is a
regular value of $f$ if $a$ is a regular value of $f_{1}$ and $b$ is a
regular value of $f_{2|f_{1}^{-1}(a)}$. 
\item[(ii)] If $f$ is a Fredholm map, then $f_{2|f_{1}^{-1}(a)}$ is
a Fredholm map for all regular values $a$ of $f_{1}$. For
a regular value $a$ of $f_{1}$, there is a Baire subset
$U\subset {\cal N}_{2}$ such that $(a,b)$ is a regular value of
$f$ for $b\in U$. 
\end{description}
}
\label{transverse3}
\end{lemma}
\vspace{.25in}
\noindent
{\bf Proof:} Let $M_{a}:=f_{1}^{-1}(a)$. We have the commuting diagram:
$$
\begin{array}{lclcl}
M_{a}&\hookrightarrow&{\cal M}&\stackrel{f_{1}}{\longrightarrow}
&{\cal N}_{1}\\
\downarrow\,f_{|M_{a}}& &\downarrow\,f& &\left|\!\right|\\
\{a\}\times{\cal N}_{2}&\hookrightarrow&{\cal N}_{1}\times{\cal
N}_{2}&\longrightarrow&{\cal N}_{1}
\end{array}
$$
which leads to the diagram of derivatives :
$$
\begin{array}{cclclclcc}
0&\rarr&T_{x}(M_{a})&\hookrightarrow &T_{x}{\cal M}&\stackrel{Df_{1}(x)}
{\longrightarrow}&T_{a}{\cal N}_{1}&\rarr &0\\
& &\downarrow\,Df_{2|M_{a}}& &\downarrow\,Df(x)&
&\left|\!\right|& &\\
0&\rarr&T_{f_{2}(x)}{\cal N}_{2}&\longrightarrow&T_{a}{\cal N}_{1}\oplus 
T_{f_{2}(x)}{\cal N}_{2}&\longrightarrow&T_{a}{\cal N}_{1}&\rarr &0
\end{array}
$$for all $x\in {\cal M}_{a}$. Clearly $Df(x)$ is surjective iff
$Df_{2}(x)$ is surjective, for all $x\in {\cal M}_{a}$. Thus
$(a, f_{2}(x))$ is a regular value of $f$ iff $f_{2}(x)$ is a regular
value for $f_{2|{\cal M}_{a}}$. Similarly, the Fredholm
statement follows because of the snake lemma, for all $x\in
{\cal M}_{a}$. \hfill $\Box$


\begin{thebibliography}{999}
\bibitem[D]{} Donaldson, S.,{\em The Seiberg-Witten Equations and
4-Manifold Topology}, Bull. A.M.S., Vol 33, No. 1, 45-71, (1996).
\bibitem[Ko]{} Kodaira, Kunihiko, {\em Harmonic Fields in Riemannian
Manifolds},Ann. of Math., 50 (1949), 587-665.
\bibitem[KM]{}Kronheimer, P., and Mrowka, T.S., {\em The Genus of
Embedded Surfaces in the Projective Plane}, Math. Res. Letters 1, 797-808
(1994).
\bibitem[PP]{}Paranjape, Kapil, and Pati V., {\em Seiberg-Witten
Invariants, an Expository Account}, (this math.dg archive)
\bibitem[W]{}Witten, E., {\em Monopoles and 4-Manifolds}, Math Res.
Letters, 1, 769-796, (1994).
\end{thebibliography}
\end{document}